\documentclass[12pt,a4paper]{article}
\makeatletter \@addtoreset{equation}{section} \makeatother
\def\N{\mbox{I\hspace{-.15em}N}}
\def\R{\mbox{I\hspace{-.15em}R}}

\newtheorem{theo}{Theorem}[section]
\newtheorem{lem}[theo]{Lemma}
\newtheorem{rem}[theo]{Remark}

\newtheorem{pro}[theo]{Proposition}
\def\u{\mathbf{u}}

\def\o{\mbox{\boldmath$\omega$\unboldmath}}

\def\vfy{\mbox{\boldmath$\varphi$\unboldmath}}

\def\O{\Omega}
\def\v{\mathbf{v}}
\def\c{\mathbf{curl}\,}
\def\p{\Vert}
\def\f{\mathbf{f}}
\def\n{\nabla}
\def\d{\mathrm{div}\,}

\def\P{\partial}
\def\F{\frac}
\def\D{\Delta}
\def\a{\alpha}

\def\t{\times}
\def\w{\mathbf{w}}
\def\z{\mathbf{z}}
\def\x{\mathbf{x}}
\def\g{\mathbf{g}}
\def\h{\mathbf{h}}
\def\r{$\rbrack$}
\def\l{$\lbrack$}
\def\hs{\hspace*{0.2cm}}
\def\mb{\mathbf}
\def\mr{\mathrm}
\def\H{\hspace*{0.7cm}} 
\def\iy{\infty}

\def\di{\displaystyle}
\def\rk{$\rbrack$}
\def\lk{$\lbrack$}
 
\def\tho{\mbox{\boldmath$\tau$\unboldmath}}
\def\ve{\varepsilon}  

\def\cs{\mr{curl}\,}
\def\tb{\textbf}
\def\G{\Gamma}
\def\ov{\overline}
\def\m{\vert}
\def\nb{\mb{n}}
\def\W{\mathcal{W}}
\begin{document}
\setcounter{page}{0}
\title{Fully non-homogeneous problem of two-dimensional second  
grade fluids.} 
\author{J. M. Bernard*}
\date{ }
\maketitle
\begin{abstract}
This article studies the solutions of a two-dimensional grade-two fluid model with a fully non-homogeneous boundary condition 
for velocity $\u$. Compared to pro-\linebreak blems with a homogeneous or tangential boundary condition, studied by many authors, we must add a boundary condition, otherwise the problem is no longer well-posed. We propose two conditions on $z=\cs(\u-\a\D\u)$, which differ according to the regularity of $z$, on the portion of 
$\P\O$ where $\a\u\,.\,\mb{n}<0$. Following the approach of V. Girault and L.R. Scott in the tangential boundary case, we split the problem into a system with a generalized Stokes problem and a transport problem. But, compared to the study of these authors, we are now led to solve  transport problems with boundary conditions. In two previous articles, we studied these transport pro-\linebreak blems. The results obtained in these articles allow us, by a fixed-point argument, to establish existence of the solutions for the fully non-homogeneous grade-two problem. Uniqueness requires the boundary condition with 
$z$ in $H^1$.
\end{abstract} 
 \vspace{8cm}
*Universit\'e d'Evry Val d'Essonne, Boulevard F. Mitterand. \\
91025 Evry Cedex, France.
 \pagebreak
\section{introduction}
\H This paper studies the stationary problem of a class of second-grade
fluids in two dimensions. The system of equations we propose to solve is:
\begin{eqnarray}
 -\nu\Delta\u+\c(\u-\alpha\Delta \u)\times \u+\nabla p=
 \f\ \mathrm{in}\ \O,\label{1eu}\end{eqnarray}
with the incompressibility condition:
\begin{equation}
 \label{1ed} \d\,\u=0\hspace{1cm} \mathrm{in}\hspace{0.3cm}\O,
\end{equation}
with adequate non-homogeneous boundary conditions.\\
\H A grade-two fluid is a non-Newtonian fluid and it is considered as an
appropriate model for the motion of a water solution of polymers,
cf. Dunn and Rajagopal [17]. The parameter $\nu$ is the viscosity and the
parameter $\a$ is a constant stress modulus, both divided by the density.
When $\a=0$, the constitutive equation reduces to that of the
Navier-Stokes equation.\\
\H The thermodynamics of fluids of grade 2 entail that $\nu$ and $\a$ be
non-negative (cf.[16]), but, since the sign of $\a$ in (\ref{1eu}) is
unimportant from a strictly mathematical point of view,
we only shall assume $\nu>0$.\\
\H Concerning fluids of grade $n$, we refer to W. Noll and C. Truesdell [27],
R.L. Fosdick and K.R. Rajagopal [18,19].\\
\H We write $\u=(u_1,u_2,0)$ in order to define the curl and the vector product.
Recall that $\c\u=(0,0,\cs\u)$, where
$$\cs\u=\F{\P u_2}{\P x_1}-\F{\P u_1}{\P x_2}.$$
\H We impose a fully non-homogeneous Dirichlet boundary condition :
\begin{equation}\label{1b1}
\u=\g\H\H \mr{on}\ \P\O\H\mr{with}\H \int_{\gamma_i}\g\,.\,\mb{n}\,ds=0,
\H 0\le i\le k,\end{equation}
  where $\gamma_i$, $0\le i\le k$,
 denotes the connected components of its boundary $\P\O$ and
 $\mb{n}=(n_1,n_2)$ denotes the unit exterior normal to the boundary $\P\O$
of $\O$.\\
\H Next, since we do not assume $\g\,.\,\nb=0$ on $\P\O$, we have to impose supplementary conditions on parts of the boundary, otherwise the problem is no longer well posed, as we shall see later.\\ 
\H Let us denote by $\G^-$ the following open portion of
$\P\O$
\begin{equation}\label{1dg-}
\G^-=\bigcup_{i\in I} \o_i,\end{equation}
where the sequence $(\o_i)_{i\in I}$ represents the set of the open
sets $\o_i$ of $\P\O$ such that 
$\a\,\g\,.\,\mb{n}<0$ almost
everywhere in $\o_i$.\\  
\H In the same way, let us denote by $\G^{0,+}$ the following open portion of $\P\O$
\begin{equation}\label{1dg0+}
\G^{0,+}=\bigcup_{j\in J} \o'_j,\end{equation} 
where the open sets $\o'_j$ of $\P\O$ are such that
$\a\,\g\,.\,\mb{n}\ge 0$ almost
everywhere in $\o'_j$. 
Let us note that these definitions imply
 $$\G^-\cap\G^{0,+}=\emptyset.$$
 \H We assume that $\G^-$ and $\G^{0,+}$ have a finite number of connected components and verify
\begin{equation}\label{1dg0g1}
\P\O=\ov{\G^-}\cup\ov{\G^{0,+}},\H
\ov{\G^-}\cap\ov{\G^{0,+}}=\{\mb{m}_1,\ldots,\mb{m}_q\},
\end{equation}
where $\mb{m}_k$, $1\le k\le q$, denotes points of the boundary $\P\O$.\\
\H First, we impose the following additional condition on $\G^-$ :
\begin{equation}\label{1b2}
(\cs(\u-\a\D\u)\u)\,.\,\mb{n}=h\H \mr{on}\ \G^{-}.\end{equation}
With this additional condition, we will obtain the existence of solutions for the fully non-homogeneous problem of two dimensional second grade fluids under rather mild assumptions on the data, but we cannot prove the uniqueness of solution. In order to obtain 
uniqueness, we are led to assume another condition on $\G^-$, which requires $\cs(\u-\a\D\u)$ in $H^1$. So, the 
boundary condition on $\G^-$ has the following simpler formulation, namely : 
\begin{equation}\label{1b3}
\cs(\u-\a\D\u)=h\H \mr{on}\ \G^{-}.
\end{equation}
With this second additional boundary condition, we will obtain existence and uniqueness for the fully non-homogeneous problem of two dimensional second grade fluids, but under stronger assumptions on the data and the boundary.
\vspace{0,3cm}\\
 \H The difficulty of this problem arises from the fact that its elliptic
term is only a Laplace operator, whereas its nonlinear term involves a
third-order derivative. Roughly, two approaches have been used to
study the grade-two problem. The first one is a method of energy estimates,
initiated by Ouazar [28] in 1981 and Cioranescu and Ouazar [11,12]. They
look for a velocity $\u$ such that $\z=\c(\u-\a\D\u)$ has $L^2$
regularity, introducing $\z$ as an auxiliary variable and discretizing
the equations of motion by Galerkin's method in the basis of the
eigenfunctions of the operator $\c\c(\u-\a\D\u)$. Cioranescu and Girault [10],
Bernard [2,3,4] and more recently Girault and
Scott [23], by using the renormalizing technique of [15], extended the results of Cioranescu and Ouazar for
both the time-dependent and steady-state grade two fluid model in two dimensions. In 2012, following the approach of Girault and Scott in [23], Bernard [6] study the steady-state grade two fluid model in convex polyhedron.\\
\H For about twenty years, authors such as Bresch and Lemoine [8,9], Costia
and Galdi [13], Galdi, Grobbelaar-Van Dalsen and Sauer [20], Galdi and
Sequeira [21], Videman [29] have used another approach: each one decomposed
the original system of equations in their own way, but all applied a
Schauder fixed point argument. In particular, Videman proves existence and
uniqueness in $W^{2,p}$ with a boundary of class $C^{1,1}$ for
sufficiently small data.\\
\H Comparing the two approaches, the method of energy estimates of
Cioranescu and Ouazar is the only one that gives existence of solutions
in two dimensions for the second grade fluids, without restriction on the
size of the data. Indeed, in the methods using a Schauder fixed point
argument, the nonlinear term is placed straight without conversions on the
right hand side and the existence of solutions is thus proven with
heavy restrictions on the size of the data and parameters.\\
\H Since the work of Cioranescu and Ouazar, the most important progresses
in the study of grade-two fluid were done by Girault and Scott in
the already quoted paper [23]. They have studied the solutions in $H^1$ of
a two-dimensional grade-two fluid model with a non-homogeneous Dirichlet
tangential boundary condition, on a Lipchitz-continuous domain, this
weak regularity of the boundary allowing for a subsequent numerical
analysis of the
model. Both for numerical purposes and solving the difficulty of a boundary
with few regularity, they developped the next variant of the method of
Cioranescu and Ouazar. The idea is to split the original system of equations
into a coupled generalized Stokes problem satisfied by $\u$ and a
transport equation satisfied by $z$. They obtained existence of solutions
without restriction on the size of the data and the constant parameters of the
fluid. A substantial part of this article was
devoted to a sharp analysis of the transport equation under weak regularity
assumptions. Uniqueness was established in a convex polygon, with adequate
restrictions on the size of the data and parameters. In addition,
they proved a difficult result that was hitherto regarded as one of the
major open questions relative to models of grade-two fluids: any solution
of the grade-two problem converges to a solution of the
Navier-Stokes equations when $\a$ tends to zero. As a result of the weak
regularity of the boundary, V. Girault and L.R. Scott proposed  
 finite element discretizations of a two dimensional grade-two fluid
 model in [24].\\
\H Another major and still open question relative to models of grade-two fluids is the
fully non-homogeneous problem. As Girault and Scott wrote in [23], this problem is not well-posed, thus implying that additional
boundary conditions should be imposed. But, as they pointed out, it was not
yet known what boundary conditions could be imposed in order to insure that
the problem is well-posed. The purpose of this paper is to give additional boundary conditions, which insures that the fully non-homogeneous grade two problem is well-posed. But, one of the main difficulty for the fully non-homogeneous problem arises from the fact that the transport equation associated is such that the normal component of the velocity does not vanish on the boundary. In this case, the transport equation has no longer a solution and a boundary condition on $\G^{-}$ is required for the transport problem to be well posed. So, the results about these transport problems obtained by Bernard in [5] for solutions in $L^2$ and [7] for the solution in $H^1$ will be a basic tool of the proofs of existence and uniqueness.\\[0.2cm]
\H After this introduction, this article is organized as follows. In section 2, each of the two initial problems $(P_I)$ and $(P_{II})$, different depending on the additional boundary condition, are split into two equivalent coupled systems consisting of a generalized Stokes problem and a transport problem known as a mixed formulation in the way of V. Girault and L.R. Scott in [23]. Section 3 is devoted to the existence of solution, obtained by a fixed point argument, of the coupled system equivalent to Problem $(P_I)$. In section 4, we extend the existence of solution to the coupled system equivalent to Problem $(P_{II})$ . Finally, in Section 5, we prove uniqueness  of the solution of Problem $(P_{II})$ in two different frameworks: a first uniqueness theorem with weak enough assumptions but with one of the conditions of uniqueness that depends on the semi norm $H^1$ of the solution $z$ and a second uniqueness theorem with more restrictive assumptions but with conditions of uniqueness only depen-\linebreak ding on the data.\\ 
  \H In order to set this problem into adequate spaces, recall some definitions of spaces and norms. For  vector-valued functions $\v=(v_1,v_2,\ldots,v_{N})$, we use special 
norms: if $1\le p\le \infty $, we set
\begin{eqnarray}\label{2vn}\p\v\p_{L^{p}(\O)^{N}}=
\p\,\m\v\m\,\p_{L^{p}(\O)},
\end{eqnarray}
where $\m\ .\ \m$ is the euclidian norm in $\R^{N}$. To simplify, we shall 
denote $\p\mathbf{v}\p_{L^{p}(\O)}$ instead 
of $\p \mathbf{v}\p_{L^{p}(\O)^{N}}$ and $\p\v\p_{W^{m,p}(\O)}$ instead of 
$\p\v\p_{W^{m,p}(\O)^N}$.\\[0.3cm]
\H We shall frequently use the scalar product of $L^2(\O)$\\
$$(f,g)=\int_{\Omega}\ f(\mathbf{x})g(\mathbf{x})\,d\mathbf{x},$$\\ 
the semi-norm of $H^1(\O)$\\
$$\m v\m_{H^1(\O)}=\p\nabla v\p_{L^2(\O)},$$\\
and the subspaces of $H^1(\O)$ and $L^2(\O)$\\
$$H_0^1(\O)=\lbrace v\in H^1(\O);v=0\hspace{0.2cm}\mathrm{on}\ \P\O\rbrace,$$
$$H(\cs,\O)=\lbrace \v\in L^2(\O)^2;\ \cs\v 
\in L^2(\O)\rbrace,$$ 
$$V=\{\v\in H^1_0(\O)^2;\ \d\v=0\}.$$
We shall often use Sobolev's imbeddings: for any real numbers $p\ge 1$,
there exists a constant $S_p$ such that
\begin{equation}\label{1dSp}
\forall v\in H^1_0(\O),\H \p v\p_{L^p(\O)}\le S_p| v|_{H^1(\O)}.
\end{equation}
When $p=2$, this reduces to Poincar\'e's inequality and $S_2$ is Poincar\'e's
constant.\\
\H For $H^1(\O)$, we recall Sobolev's imbeddings:
\begin{equation}\label{1dSp*}
\forall v\in H^1(\O),\H \p v\p_{L^p(\O)}\le {S^*}_p\p v\p_{H^1(\O)}.
\end{equation}
 \H Let $\Gamma'$ be an open part of the boundary $\P\O$ of class $C^{0,1}$
and, for $r>2$, $T_{1,r}^{\Gamma'}$ be the mapping $v\mapsto v_{|\Gamma'}$
defined on $W^{1,r}(\O)$.
We denote by $W^{1-\F{1}{r},r}(\Gamma')$ (see [26]) the space
$T_{1,r}^{\Gamma'}(W^{1,r}(\O))$ which is equipped with the norm:
\begin{equation}\label{1dn1/2}
\p\varphi\p_{W^{1-1/r,r}(\Gamma')}=
\inf\{\p v\p_{W^{1,r}(\O)},\ v\in W^{1,r}(\O)\ \mr{and}\ v_{|\Gamma'}=\varphi\}.
\end{equation}
\H For fixed $\u$ in $H^1(\O)^2$, let us introduce the space 
 \begin{equation}\label{1dXu}
 X_{\u}(\O)=\{z\in L^2(\O),\ \u\,.\,\n z\in L^2(\O)\},\end{equation}
 which is a Hilbert space equipped with the norm
 \begin{equation}\label{1dnXu}
 \p z\p_{\u}=(\p z\p_{L^2(\O)}^2+ \p\u\,.\,\n z\p_{L^2(\O)}^2)^{1/2}.
 \end{equation}
  In the same way we define
  $$Y_{\u}(\O)==\{z\in L^2(\O),\ \u\,.\,\n z\in L^1(\O)\}.$$
  We recall a theorem ( see [5]) concerning the normal component of boundary values of $(z\u)$
where $z$ belongs to $Y_{\u}(\O)$.
\begin{theo}\label{1tg'n}
Let $\O$ be a Lipschitz-continuous domain of $\R^d$, let $\u$ belong to
$H^1(\O)^d$ with $\d\u=0$ in
$\O$ and let $r>d$ be a real number. We denote by $r'$ the
real number defined by: $\di\F{1}{r}+\di\F{1}{r'}=1$. The mapping
$\gamma'_{\mb{n}}:\ z\mapsto (z\u)\,.\,\mb{n}_{|\P\O}$ defined on
$\mathcal{D}(\overline{\O})^d$ can be extended by continuity to a linear and continuous
mapping, still denoted by $\gamma'_{\mb{n}}$, from $Y_{\u}(\O)$ into
$W^{-1/r',r'}(\P\O)$.\end{theo}
\H From this theorem and with a density argument, we derive the following
  Green's formula: let $r>d$ be a real number and let $\u$ be in
  $H^1(\O)^d$ with $\d\u=0$ in $\O$,
      \begin{equation} \label{1green}
  \forall z\in Y_{\u}(\O),\ \forall\varphi\in W^{1,r}(\O),
  \ \int_{\O}z(\u\,.\,\n\varphi)\,d\mb{x}
  +\int_{\O}\varphi(\u\,.\,\n z)\,d\mb{x}
  =<(z\u)\,.\,\mb{n},\varphi>_{\P\O}.\end{equation}
  \H Let $\G_0$ and $\G_1$ be two non-empty open parts of $\P\O$ that have a finite number of connected components and verify 
  $$\G_0\cap\G_1=\emptyset,\H
\P\O=\ov{\G_0}\cup\ov{\G_{1}},\H
\ov{\G_0}\cap\ov{\G_{1}}=\{\mb{m}_1,\ldots,\mb{m}_q\}.$$
\H We introduce
the space $W^{-1/r',r'}(\G_0)=(W^{1-1/r,r}_{00}(\G_0))'$, where 
\begin{equation}\label{1dW00r}
W^{1-1/r,r}_{00}(\G_0)=\{v_{|\G_0},\ v\in W^{1,r}(\O),\
v_{|\G_1}=0\},\end{equation}
and we denote $<\,.\,,\,.\,>_{\G_0}$ the duality pairing between these
two spaces. Note that if $z\in Y_{u}(\O)$, then
$(z\u)\,.\,\mb{n}_{|\G_0}\in W^{-1/r',r'}(\G_0)$ and, in the same way
as previously, we have the Green's formula : $\forall z\in Y_{\u}(\O)$, 
$\forall\vfy\in W^{1,r}(\O)$, with $\vfy_{|\G_1}=0$, 
$\forall\u\in H^1(\O)^d$ with $\d\,\u=0$ in $\O$,
\begin{equation}\label{1green2}
 \int_{\O}z(\u\,.\,\n\varphi)\,d\mb{x}
  +\int_{\O}\varphi(\u\,.\,\n z)\,d\mb{x}
  =<(z\u)\,.\,\mb{n},\varphi>_{\G_0}.\end{equation}
Then, we can define the following space : 
\begin{equation}\label{1dXug0}
X_{\u}(\G_0)=\{ z\in X_{\u},\ (z\u)\,.\,\mb{n}_{|\G_0}=0\}.\end{equation}
Finally, we recall a basic result of [5]. We apply this result in the particular case where $d=2$ and therefore, for $1\le k\le q$,  the sets $K_k$ are points $\mb{m_k}$ of the boundary.
\begin{pro}\label{2pgreenfg-} Let $\O$ be a Lipschitz-continuous domain
  of $\R^2$, let $\u$ be in $H^1(\O)^d$ with $\d\u=0$ in $\O$ and let
  $\G^-$ and $\G^{0,+}$ be defined by (\ref{1dg-}) and (\ref{1dg0+}),
  verifying (\ref{1dg0g1}). 
Let $z$ belong to
 $X_{\u}(\G^-)$ and $w$ to $X_{\u}(\G^{0,+})$ . Then, $z$ and $w$ verify the following inequalities
\begin{equation}\label{2greenfg-+} 
\int_{\O}\a\u\,.\,\n z)\,z\,d\x\ge
0,\H\int_{\O}(\a\u\,.\,\n w)\,w\,d\x\le 0.\end{equation}\end{pro}
\H Finally, we introduce the spaces
\begin{equation}\label{1dW}
W=\{\v\in (H^1(\O))^2;\ \mr{curl}\D\v\in L^2(\O)\},\end{equation}
\begin{equation}\label{1dW1}
W_1=\{\v\in (W^{1,\iy}(\O))^2;\ \mr{curl}\D\v\in H^1(\O)\},\end{equation}
in which we shall look for the velocity $\u$. According to the two additional boundary conditions on $\G^-$, we shall study two different problems. The first one, called $(\mr{P_I})$, with condition (\ref{1b2}) : \\ 
\H find $(\u,p)\in W\t L^2_0(\O)$ such that
$$(\mr{P_I})\left\lbrace\begin{array}{cc}
-\nu\D\u+\c(\u-\a\D \u)\t \u+\n p=
 \f\ &\mathrm{in}\ \O, \\
  \d\,\u=0\hspace{1cm} &\mathrm{in}\hspace{0.3cm}\O,\\
  \u=\g\H\H &\mr{on}\ \P\O,\\
  (\cs(\u-\a\D\u)\u)\,.\,\mb{n}=h\H& \mr{on}\ \G^{-}.  
 \end{array}\right.$$
 \H The second one, called $(\mr{P_{II}})$, with condition (\ref{1b3}) : \\ 
\H find $(\u,p)\in W_1\t L^2_0(\O)$ such that
$$(\mr{P_{II}})\left\lbrace\begin{array}{cc}
-\nu\D\u+\c(\u-\a\D \u)\t \u+\n p=
 \f\ &\mathrm{in}\ \O, \\
  \d\,\u=0\hspace{1cm} &\mathrm{in}\hspace{0.3cm}\O,\\
  \u=\g\H\H &\mr{on}\ \P\O,\\
  \cs(\u-\a\D\u)=h\H& \mr{on}\ \G^{-}.  
 \end{array}\right.$$
 \H For each problem, we shall make adequate assumptions on the data to define an equi-\linebreak valent formulation.
\section{Equivalent formulations for Problems $\mr{(P_I)}$ and $\mr{(P_{II})}$} 
\H Following the approach of [23], we shall establish a mixed formulation
of the two problems.
In this subsection, the assumptions on the data are: $\O$ is a bounded domain in
$\R^2$, with lipschitz-continuous boundary $\P\O$, $\f$ is a given
function in $H(\cs;\O)$, $\g$ is a given vector field in
$H^{1/2}(\P\O)^2$  such that
$\int_{\gamma_i}\g\,.\,\mb{n}\,ds=0$ for $0\le i\le k$, where $\gamma_i$,
$0\le i\le k$, denotes the connected components of the boundary $\P\O$ of
$\O$, $h$ is a given function in $W^{-1/r',r'}(\G^{-})$, where the real number $r'$ is defined by $\F{1}{r}+\F{1}{r'} =1$ 
from a real number $r>2$,
and $\nu>0$ and $\a$ are
two given real constants.\\
\H Let $(\u,p)\in W\t L^2_0(\O)$ be a solution of
$\mr{(P_I)}$ and introduce the auxiliary variables:
\begin{equation}\label{1dz}
z=\mr{curl}(\u-\a\D\u),\H \z=(0,0,z).\end{equation}
Note that $z\in L^2(\O)$ and
\begin{equation}\label{1divz}
\d\, \z=0.\end{equation}
With these notations, we write
(\ref{1eu}) as:
\begin{equation}\label{1esg}
 -\nu\D\u+\z\t\u+\n p=\f\H\mr{in}\ \O,\end{equation}
 that is with (\ref{1ed}) and (\ref{1b1}) a non-standard generalized Stokes equation.
 Taking the curl of (\ref{1esg}) in the sense of distributions, we
 obtain $$-\nu\,\mr{curl}\D\u+\u\,.\,\n z=\mr{curl}\,\f,$$
 that we can write as a transport equation verified by z:
 \begin{equation}\label{1et}
 \nu z +\a \u\,.\,\n z=\nu\,\cs\u+\a\,\cs\f.\end{equation}
 Moreover, the boundary condition (\ref{1b2}) can be written
 \begin{equation}\label{1b2'}(z\u)\,.\,\mb{n}=h\H\H\mr{on}\H \G^-.\end{equation}
 \begin{rem} Assume that $\u$ belongs to $W$ and
 that $\f$ belongs to $H(\c;\O)$. If $z$ is a
 solution of the transport  equation (\ref{1et}) in $L^2(\O)$,
 then $\u\,.\,\n z$ belongs to $L^2(\O)$. Finally,
 we obtain that $z$ is in $X_{\u}$. Under these assumptions,
 $(z\u)\,.\,\mb{n}=(\cs(\u-\a\D\u)\u)\,.\,\mb{n}$ belongs to
  $W^{-1/r',r'}(\G^{-})$.
 \end{rem}
 \H Conversely, let $(\u,p,z)\in (H^1(\O))^2
 \t L^2_0(\O)\t L^2(\O)$ be a
 solution of (\ref{1esg}), (\ref{1ed}), (\ref{1b1}), (\ref{1et}) and (\ref{1b2'}) and $\z=(0,0,z)$.
 Then $\z$ satisfies (\ref{1divz}) and taking the curl of (\ref{1esg}) in the
 sense of distributions yields:
 $$-\nu\,\mr{curl}\D\u+\u\,.\,\n z=\mr{curl}\,\f.$$
 Next, multiplying by $\a$ and comparing with (\ref{1et}), we obtain:
 $$z=\cs(\u-\a\D\u).$$
 Therefore $\u$ belongs to $W$ and substituting the expression of $z$
 into (\ref{1esg}) shows that $(\u,p)$ is a solution of the
 original equations (\ref{1eu})-(\ref{1b2}). This is summarised in the
 following lemma.
 \begin{lem}\label{1leqf}  Problem $\mr{(P_I)}$
 with $(\u,p)$ in $W\t L^2_0(\O)$ is
 equivalent to: Find $(\u,p,z)$ in $(H^1(\O) )^2 \t L^2_0(\O)\t L^2(\O)$
 solution of the generalized Stokes problem (\ref{1esg}), (\ref{1ed}),
 (\ref{1b1}) and the transport problem (\ref{1et}), (\ref{1b2'}), namely:
\begin{eqnarray}
 -\nu\D\u+\z\t\u+\n p=\f\H\mr{in}\ \O,\nonumber\\
  \d\,\u=0\hspace{1cm} \mathrm{in}\hspace{0.3cm}\O,\nonumber\\
\u=\g\hspace{0.8cm} \mr{on}\ \P\O,\label{1eqe}\\
  \nu z +\a \u\,.\,\n z=\nu\,\cs\u+\a\,\cs\f\H\mr{in}\ \O,\nonumber\\
 (z\u)\,.\,\mb{n}=h\H\mr{on}\ \G^-.\nonumber
 \end{eqnarray}
 \end{lem}
 \H In the same way, we establish an equivalent formulation for Problem $\mr{(P_{II})}$, with the following boundary condition on $\G^-$ :
 \begin{equation}\label{1b3'}
z=h\H \mr{on}\ \G^{-}.
\end{equation}
 \begin{lem}\label{2leqf}  Problem $\mr{(P_{II})}$
 with $(\u,p)$ in $W_1\t L^2_0(\O)$ is
 equivalent to: Find $(\u,p,z)$ in $(W^{1,\iy}(\O) )^2 \t L^2_0(\O)\t H^1(\O)$
 solution of the generalized Stokes problem (\ref{1esg}), (\ref{1ed}),
 (\ref{1b1}) and the transport problem (\ref{1et}), (\ref{1b3'}), namely:
\begin{eqnarray}
 -\nu\D\u+\z\t\u+\n p=\f\H\mr{in}\ \O,\nonumber\\
  \d\,\u=0\hspace{1cm} \mathrm{in}\hspace{0.3cm}\O,\nonumber\\
\u=\g\hspace{0.8cm} \mr{on}\ \P\O,\label{2eqe}\\
  \nu z +\a \u\,.\,\n z=\nu\,\cs\u+\a\,\cs\f\H\mr{in}\ \O,\nonumber\\
 z=h\H\mr{on}\ \G^-.\nonumber
 \end{eqnarray}
 \end{lem}
 \section{Existence of a solution for Problem (\ref{1eqe})}
\subsection{Estimates for the solution of the generalized Stokes problem}
\H In this subsection, the assumptions on the data are the same as the
previous section 2. For a given $\z$ in $L^2(\O)^3$, the generalized Stokes problem
(\ref{1esg}), (\ref{1ed}), (\ref{1b1}) has the following variational
formulation: Find $(\u(\z),p(\z))$ in $H^1(\O)^2\t L^2_0(\O)$, such that
\begin{eqnarray}
\forall\v\in H^1_0(\O)^2,\H a_{\z}(\u(\z),\v)+b(\v,p(\z))=(\f,\v),
\label{2sg1}\\
\forall q\in L^2_0(\O),\H b(\u(\z),q)=0,\label{2sg2}\hspace*{4cm}\\
\u(\z)=\g\hspace{0.5cm} \mr{on}\ \P\O\H\H\mr{with}\H\int_{\gamma_i}\g\,.\,\mb{n}\,ds=0,\ 0\le i\le k,
\label{2sg3}
\end{eqnarray}
where$$ a_{\z}(\w,\v)=\nu(\n\w,\n\v)+(\z\t\w,\v),$$
$$b(\v,q) =-(q,\d\v).$$
In the same way as in [23], we define two liftings of $\g$: first, we
lift $\g$ by $\w_{\g}$ solution in $H^1(\O)^2$ of the non-homogeneous
Stokes problem:
\begin{equation}\label{2lwg}
-\D\w_{\g}+\n p_{\g}=\mb{0}\H\mr{and}\H\d\w_{\g}=0\H\mr{in}\ \O,\H
 \w_{\g}=\g\H\mr{on}\ \P\O.\end{equation}
 Under the assumption: $\int_{\gamma_i}\g\,.\,\mb{n}\,ds=0$ for $0\le i\le k$, this
 problem has a unique solution and it satisfies the bound (see [22]):
 \begin{equation}
 \label{2bwg}
 \p\w_{\g}\p_{H^1(\O)}\le T\p\g\p_{H^{1/2}(\P\O)}.\end{equation}
 \H To show the existence of solutions without restriction on the data,
 we need to construct an adequate lifting $\u_{\g}$ (see [23]) in the same way as
 a lemma by Leray and Hopf in the case of the nonhomogeneous Navier-Stokes
 equations.
 \begin{theo}\label{2lug} Let $\O$ be a lipschitz-continuous domain and let $\gamma_i$,
 $0\le i\le k$, denote the connected components of its boundary $\P\O$.
 There exists a continuous non-increasing function $L:\R^+\mapsto \R^+$,
 that tends to infinity as its argument tends to zero, such that for any
 real number $\varepsilon>0$ and for all function $\g$ in $H^{1/2}(\P\O)^2$
 satisfying
 $$\int_{\gamma_i}\g\,.\,\mb{n}\,ds=0,\H0\le i\le k,$$
 there exists a lifting function $\u_{\g}$ in $H^1(\O)^2$ with:
 $$\d \u_{\g}=0\H\mr{in}\ \O,\H\H\u_{\g}=\g\H\mr{on}\ \P\O,$$
 \begin{equation}\label{2mug1}
 \p\u_{\g}\p_{H^1(\O)}\le L\left(\F{\varepsilon}{\p \g\p_{H^{1/2}(\P\O)}}
 \right)\p \g\p_{H^{1/2}(\P\O)},\end{equation}
  \begin{equation}\label{2mug2}
  \forall\v\in H^1_0(\O)^2,\H\p |\u_{\g}|\,|\v|\p_{L^2(\O)}\le
  \varepsilon |\v|_{H^1(\O)}.\end{equation}
  \end{theo}
  These liftings allow us to show the following lemma.
  \begin{lem}\label{2luzb} Let $\O$ be Lipchitz-continuous, $\nu>0$, $\f\in L^2(\O)^2$
  and $\g\in H^{1/2}(\P\O)^2$ satisfying the second part of (\ref{1b1}).
  For any $\z$ in $L^2(\O)^3$, the generalized Stokes problem
  (\ref{2sg1})-(\ref{2sg3}) has a unique solution
  $(\u(\z),p(\z))$ in $H^1(\O)^2\t L^2_0(\O)$. This solution
  satisfies the following bounds:
  \begin{equation}\label{2muz1}
   \p\u(\z)\p_{H^1(\O)}\le \F{S_2\sqrt{S_2^2+1}}{\nu}\p\f\p_{L^2(\O)}
   +T\p\g\p_{H^{1/2}(\P\O)}\left(1+
   \F{S_4S^*_4\sqrt{S_2^2+1}}{\nu}\,
   \p\z\p_{L^2(\O)}\right);\end{equation}
  \begin{eqnarray}
  \forall\varepsilon>0,\H\p\u(\z)\p_{H^1(\O)}\le
  \F{S_2\sqrt{S_2^2+1}}{\nu}\p\f\p_{L^2(\O)}\hspace{4cm} \nonumber\\
  +(1+ \sqrt{S_2^2+1}\,)L\left(\F{\varepsilon}{\p \g\p_{H^{1/2}(\P\O)}}
 \right)\p \g\p_{H^{1/2}(\P\O)}
 + \F{\sqrt{S_2^2+1}}{\nu}\,\varepsilon
 \p\z\p_{L^2(\O)},\label{2muz2}
 \end{eqnarray}
 \begin{equation}\label{2mpz}
 \p p(\z)\p_{L^2(\O)}\le\F{1}{\beta}(S_2\p\f\p_{L^2(\O)}+\nu T
 \p \g\p_{H^{1/2}(\P\O)}+S_4S^*_4\p\u(\z)\p_{H^1(\O)}
 \p\z\p_{L^2(\O)}),\end{equation}
  where $\beta>0$ is the isomorphism constant of the divergence operator,
  as given in formula (\ref{2ibeta}) below, $S_p$ and $S^*_p$ are
  defined in
  (\ref{1dSp}) and (\ref{1dSp*}) respectively and $T$ is defined in (\ref{2bwg}).
  \end{lem}
 \tb{Proof.} Let $\u^*_g$ be any lifting of $\g$ such that
 $\u_0=\u-\u_{\g}^*$ belongs to $V$. Then (\ref{2sg1})-(\ref{2sg3}) is
 equivalent to: Find $\u_0\in V$ such that:
 \begin{equation}\label{2evu0}
 \forall \v\in V,\H a_{\z}(\u_0,\v)=(\f,\v)- a_{\z}(\u_{\g}^*,\v).
 \end{equation}
 For fixed $\z$ in $L^2(\O)^3$, the bilinear form $a_{\z}$ is elliptic on
 $H^1_0(\O)^2\t H^1_0(\O)^2$ since $(\z\t\v,\v)=0$, and it is continuous on
 $H^1(\O)^2\t H^1(\O)^2$ since
 $$\forall\u,\ \v\in H^1(\O)^2,\H |(\z\t\u,\v)|\le \p\z\p_{L^2(\O)}
 \p\u\p_{L^4(\O)}\p\v\p_{L^4(\O)}.$$
 Moreover, the linear form $\v\mapsto (\f,\v)- a_{\z}(\u_{\g}^*,\v)$ is
 continuous on $H^1_0(\O)^2$. Therefore, (\ref{2evu0}) has a unique solution
 $\u_0\in V$ and in turn this implies that (\ref{2sg1})-(\ref{2sg3}) has
 a unique solution $(\u(\z),p(\z))$ in $H^1(\O)^2\t L^2_0(\O)$.\\
 \H Taking for $\u_{\g}^*$ the lifting $\w_{\g}$ defined by (\ref{2lwg}),
 the choice $\v=\u_0$ in (\ref{2evu0}) yields
 $$|\u_0|_{H^1(\O)}\le \F{1}{\nu}(S_2\p\f\p_{L^2(\O)}+S_4S^*_4
 \p\z\p_{L^2(\O)}\p\w_{\g}\p_{H^1(\O)}).$$
 Then using Poincar\'e's
 constant, the triangle inequality and (\ref{2bwg}), we obtain
 (\ref{2muz1}). For the second bound, we take for $\u_{\g}^*$ the lifting
 $\u_{\g}$ of  Theorem \ref{2lug}. Then choosing again $\v=\u_0$ in
 (\ref{2evu0}) and using (\ref{2mug2}), we derive
 $$|\u_0|_{H^1(\O)}\le \F{S_2}{\nu}\p\f\p_{L^2(\O)}+|\u_{\g}|_{H^1(\O)}+
 \F{\varepsilon}{\nu}\p\z\p_{L^2(\O)}.$$
 Then using again Poincar\'e's
 constant and the triangle inequality  and owing to
 (\ref{2mug1}), we obtain (\ref{2muz2}).\\
 \H Concerning $p(\z)$, it follows from the isomorphism properties of the
 divergence (cf for instance [22]) that there exists a unique $\v_p$ in
 $H^1_0(\O)^2$ such that
 $$\d \v_p=p(\z)\H\H \mr{in}\ \O,$$
 \begin{equation}\label{2ibeta}
 |\v_p|_{H^1(\O)}\le\F{1}{\beta}\p p(\z)\p_{L^2(\O)},\end{equation}
 $$\forall\w\in V,\H\H(\n\v_p,\n\w)=0.$$
 Then taking $\v_p$ for test function in (\ref{2sg1}), since
 $(\n\u(\z),\n\v_p)=(\n\w_{\g},\n\v_p)$, we obtain:
 $$\p p(\z)\p_{L^2(\O)}^2=(\z\t\u(\z),\v_p)+\nu(\n\w_{\g},\n\v_p)-
 (\f,\v_p).$$
 Therefore, applying (\ref{2ibeta}), we derive (\ref{2mpz}).
 \hfill$\diamondsuit$ 
 \subsection{Existence of a solution for Problem (\ref{1eqe})}
 \H In this subsection, we need additional assumptions on $h$ to insure
the existence of a lifting of $h$. We assume that
\begin{equation}\label{2hh}
h\in L^1(\G^-)\H\mr{and}\H
(\di\F{h}{\g\,.\,\mb{n}})_{|\G^-}\in W^{1-1/t,t}(\G^-)\H\mr{for}\ t>2.
\end{equation}
Since $(\di\F{h}{\g\,.\,\mb{n}})_{|\G^-}$ belongs to
$W^{1-\F{1}{t},t}(\G^-)$, there
exists a function $z_h\in W^{1,t}(\O)$ such that $z_h=
\di\F{h}{\g\,.\,\mb{n}}$ on $\G^-$ with
\begin{equation}\label{2nzh}
\p z_h\p_{W^{1,t}(\O)}=
\p \di\F{h}{\g\,.\,\mb{n}}\p_{W^{1-\F{1}{t},t}(\G^-)}
\end{equation}
Since $z_h$, $h$, $\mb{n}$ and $\g=\u_{|\P\O}$ are defined almost everywhere on $\G^-$, we obtain
\begin{equation}\label{2zh}
(z_h\u)\,.\,\mb{n}=h\H\H\mr{on}\H\G^-.\end{equation}
The following existence theorem is a basic result of the article.
\begin{theo}\label{2te}
   Let $\O$ be lipschitz-continuous. For all real numbers $\nu $, $\a$
   and $t$ with $\nu>0$ and $t>2$, all $\f\in H(\cs;\O)$,
   all $\g\in H^{1/2}(\P\O)^2$ satisfying the second part of (\ref{1b1}), such that $\G^-$ and $\G^{0,+}$ defined by (\ref{1dg-}) and (\ref{1dg0+}) 
   verify (\ref{1dg0g1}), and  all $h\in L^1(\G^-)$ verifying (\ref{2hh}), there exists at least one solution
   $(\u,p,z)$ for Problem (\ref{1eqe}) and
   this solution satisfies the following estimates: 
   \begin{equation}
   \p z\p_{L^2(\O)}\le \F{2|\a|}{\nu}\p\cs\f\p_{L^2(\O)}
   +4\p z_h \p_{L^2(\O)}+C(\a,\nu,\f,\g,h)
   \label{2mz},\end{equation}
  \begin{equation}\label{2mu}
   \p\u\p_{H^1(\O)}\le \F{S_2\sqrt{S_2^2+1}}{\nu}\p\f\p_{L^2(\O)}
   +T\p\g\p_{H^{1/2}(\P\O)}\left(1+
   \F{S_4S^*_4\sqrt{S_2^2+1}}{\nu}\,
   \p\z\p_{L^2(\O)}\right),\end{equation}
 \begin{equation}\label{2mp}
 \p p\p_{L^2(\O)}\le\F{1}{\beta}(S_2\p\f\p_{L^2(\O)}+\nu T
 \p \g\p_{H^{1/2}(\P\O)}+S_4S^*_4\p\u\p_{H^1(\O)}
 \p\z\p_{L^2(\O)}),\end{equation}
  where $z_h$, $C(\a,\nu,\f,\g,h)$, $T$, $\beta$, $S_p$ and
  $S^*_p$ are
  defined in (\ref{2nzh}), (\ref{3dC}), (\ref{2bwg}), (\ref{2ibeta}), (\ref{1dSp}) and (\ref{1dSp*})
  respectively.
  \end{theo}         
\tb{Proof.}  Let us define a sequence $(z^*_n)$ of functions $z^*_n\in X_{\u}(\G^-)$, $n\in\N$, by recurrence, where $X_{\u}(\G^-)$
is defined by (\ref{1dXug0}).
We set $z_0^*=0$ and assume that the function $z^*_n\in X_{\u}(\G^-)$ is given for $n\in \N$. First, setting
\begin{equation}\label{3dzn} 
z_n=z^*_n+z_h,\end{equation}
 we denote $(\u(z_n),p(z_n))$ the unique solution in $H^1(\O)^2\t L^2_0(\O)$ of the generalized Stokes problem 
\begin{equation}\label{3eqe}\left\lbrace\begin{array}{c}
 -\nu\D\u+\z_n\t\u+\n p=\f\H\mr{in}\ \O,\\
  \d\,\u=0\hspace{1cm}\mathrm{in}\hspace{0.3cm}\O,\\
\u=\g\H\H \mr{on}\ \P\O.\end{array}\right.\end{equation}
Second, we define $z_{n+1}^*\in X_{\u}(\G^-)$ 
as the unique solution of the transport problem (see Theorem 3.3 in [5])
\begin{equation}\label{3eqt}\left\lbrace\begin{array}{c}
 \nu z^*_{n+1} +\a \u(z_n)\,.\,\n z_{n+1}^*=
 \nu\,\cs\u(z_n)+\a\,\cs\f-\nu z_h-
 \a \u(z_n)\,.\,\n z_h \H\mr{in}\ \O,\\
 (z_{n+1}^*\u(z_n))\,.\,\mb{n}=0\H\H\mr{on}\ \P\O^-.\end{array}\right. 
 \end{equation}
 \H Since $z_{n+1}^*$ belongs to $X_{\u}(\G^-)$, the basic result of Proposition \ref{2pgreenfg-} implies
 $$\int_{\O}(\a \u(z_n)\,.\,\n z_{n+1}^*)z_{n+1}^*\,d\x\ge 0.$$ 
 Then, taking the scalar product of both sides of the first equation of (\ref{3eqt}) with $z_{n+1}^*$ 
 yields 
 $$\nu\p z_{n+1}^*\p^2_{L^2(\O)}\le \nu\,(\cs\u(z_n),z_{n+1}^*)+ \a\,(\cs\f,z_{n+1}^*)-\nu(z_h,z_{n+1}^*)-\a(\u(z_n)\,.\,\n z_h,z_{n+1}^*).$$
 Hence, we derive
 $$\p z_{n+1}^*\p_{L^2(\O)}\le \p\cs\u(z_n)\p_{L^2(\O)}+ \F{|\a|}{\nu}\,\p\cs\f\p_{L^2(\O)}+\p z_h\p_{L^2(\O)}+\F{|\a|}{\nu}\,\p\u(z_n)\,.\,\n z_h\p_{L^2(\O)}.$$
 Since $z_h$ belongs to $W^{1,t}(\O)$ and owing to (\ref{1dSp*}), we derive
 $$\p\u(z_n)\,.\,\n z_h\p_{L^2(\O)}\le S^*_{\F{2t}{t-2}}\p z_h\p_{W^{1,t}(\O)}\p\u(z_n)\p_{H^1(\O)}.$$
 Substituting this bound yields
 $$\p z_{n+1}^*\p_{L^2(\O)}\le (\sqrt{2}+ \F{|\a|}{\nu}S^*_{\F{2t}{t-2}}\p z_h\p_{W^{1,t}(\O)})\p\u(z_n)\p_{H^1(\O)}+ \F{|\a|}{\nu}\,\p\cs\f\p_{L^2(\O)}
 +\p z_h\p_{L^2(\O)}.$$
 Next, using the basic bound of $\p\u(z_n)\p_{H^1(\O)}$ given by (\ref{2muz2}), considering (\ref{3dzn}) and 
 setting $K=\sqrt{2}+ \F{|\a|}{\nu}S^*_{\F{2t}{t-2}}\p z_h\p_{W^{1,t}(\O)}$ , we obtain
 \begin{eqnarray*}
 \p z_{n+1}^*\p_{L^2(\O)}\le K
 (\F{S_2\sqrt{S_2^2+1}}{\nu}\p\f\p_{L^2(\O)}
  +(1+ \sqrt{S_2^2+1}\,)L\left(\F{\varepsilon}{\p \g\p_{H^{1/2}(\P\O)}}
 \right)\p \g\p_{H^{1/2}(\P\O)})\\
 +\F{|\a|}{\nu}\,\p\cs\f\p_{L^2(\O)}+(1+K\F{\sqrt{S_2^2+1}}{\nu}\,\ve)\p z_h\p_{L^2(\O)}+K\F{\sqrt{S_2^2+1}}{\nu}\,\ve\p z_n^*\p_{L^2(\O)}.
 \end{eqnarray*}
 We choose
 $$\ve=\F{\nu}{2\sqrt{S_2^2+1}\,K}=\F{\nu}{2\sqrt{S_2^2+1}(\sqrt{2}+ \F{|\a|}{\nu}S^*_{\F{2t}{t-2}}\p z_h\p_{W^{1,t}(\O)})}$$
 and we set
 \begin{eqnarray} 
   C(\a,\nu,\f,\g,h)
   =(2\sqrt{2}+2\F{|\a|}{\nu} S^*_{\F{2t}{t-2}}\p z_h\p_{W^{1,t}(\O)})
      \left(\F{S_2\sqrt{S_2^2+1}}{\nu}\p\f\p_{L^2(\O)}\right.\nonumber\\
      \left. +(1+\sqrt{S_2^2+1}\,) L(\F{\nu^2}{2\sqrt{S_2^2+1}(\sqrt{2}\nu+|\a| S^*_{\F{2t}{t-2}}
   \p z_h\p_{W^{1,t}(\O)})\p\g\p_{H^{1/2}(\P\O)}}
   )\p\g\p_{H^{1/2}(\P\O)}\right),\H\label{3dC}\end{eqnarray}
   where $L(.)$ is defined in Theorem \ref{2lug}.
   Then, we obtain the following inequality
   \begin{eqnarray*}\p z_{n+1}^*\p_{L^2(\O)}\le \F{1}{2}\p z_{n}^*\p_{L^2(\O)}+\F{|\a|}{\nu}\,\p\cs\f\p_{L^2(\O)}+\F{3}{2}\p z_h\p_{L^2(\O)}
   +\F{1}{2}C(\a,\nu,\f,\g,h),\end{eqnarray*}
   which implies, by a recurrence argument, that $z_n^*$ is uniformly bounded in $L^2(\O)$
   \begin{equation}\label{3bzn*} 
   \forall n\in\N,\ \p z_{n}^*\p_{L^2(\O)}\le 2\F{|\a|}{\nu}\,\p\cs\f\p_{L^2(\O)}+3\p z_h\p_{L^2(\O)}+C(\a,\nu,\f,\g,h).\end{equation}
   Hence, from (\ref{3dzn}), (\ref{2muz1}) and (\ref{2mpz}), we derive that $u(z_n)$ and $p(z_n)$ are uniformly bounded in $H^1(\O)^2$ and $L^2(\O)$, 
   respectively. Moreover, considering that
   $$\a\u(z_n)\,.\,\n z^*_{n+1}=\nu\,\cs\u(z_n)+\a\,\cs\f-\nu z_h-
 \a \u(z_n)\,.\,\n z_h -z_{n+1}^*$$
 and the bound
 $$\p\u(z_n)\,.\,\n z_h\p_{L^2(\O)}\le S^*_{\F{2t}{t_2}}\p\u(z_n)\p_{H^1(\O)}\p z_h\p_{W^{1,t}(\O)},$$
 we obtain that $\u(z_n)\,.\,\n z^*_{n+1}$ is uniformly bounded in $L^2(\O)$. Therefore, there  
 exists a subsequence, still denoted  by the index $m$, and four
  functions $z^*\in L^2(\O)$, $\u\in H^1(\O)^2$, $p\in L^2(\O)$, $l\in L^2(\O)$
  such that
  $$\lim_{n\to\iy}z_n^*=z^*\H\mr{weakly\ in}\ L^2(\O),$$
  $$\lim_{n\to\iy}\u(\z_n)=\u  \H\mr{weakly\ in}\ H^1(\O)^2,$$
  $$\lim_{n\to\iy}p(\z_n)=p \H\mr{weakly\ in}\ L^2(\O),$$
  $$\lim_{n\to\iy}\u(z_n)\,.\,\n z^*_{n+1}=l \H\mr{weakly\ in}\ L^2(\O).$$
  \H The weak convergence of $\u(z_n)$ in $H^1(\O)^2$ implies that for all
  real $p<\iy$
  $$\lim_{n\to\iy}\u(\z_n)=\u  \H\mr{in}\ L^p(\O)^2.$$
  Hence, we derive that $\u(z_n)\,.\,\n z^*_{n+1}$ converge to $\u\,.\,\n z^*$ in $\mathcal{D}(\O)'$, which gives 
  $l=\u\,.\,\n z^*$. Setting $z=z^*+z_h$, we obtain that $z_n$ converge to $z$ in $L^2(\O)$ weakly. These convergences allow us to pass to the limit in the generalized Stokes problem (\ref{3eqe}) and in the transport equation (\ref{3eqt}). Thus $(\u,p)$ is a solution in $H^1(\O)^2\t L^2_0(\O)$ of the generalized 
 Stokes problem (\ref{1esg}), (\ref{1ed}, (\ref{1b1}) and $z$ is a solution in $L^2(\O)$ of the transport equation (\ref{1et}).\\
 \H It remains to prove the condition (\ref{1b2}). From the Green's formula (\ref{1green2}) with $\G_0=\G^-$ and $$(z_{n+1}^*\u(z_n))\,.\,\mb{n}_{|\G^-}=0,$$ we derive
 $\forall\vfy\in W^{1,r}(\O)$, with $\vfy_{|\G^{0,+}}=0$,  
 $$(z_{n+1}^*\u(z_n),\n\vfy)+(\vfy\u(z_{n}),\n z_{n+1}^*)= <(z^*_{n+1}\u(z_n))\,.\,\mb{n},\varphi>_{\G^-}=0.$$
 Using the above convergence, we can pass to the limit and we obtain
 $$\forall\vfy\in W^{1,r}(\O),\ \mr{with}\ \vfy_{|\G^{0,+}}=0,\ (z^*\u,\n\vfy)+(\vfy\u,\n z^*)=0,$$
 which implies, again with Green's formula (\ref{1green2}), $<(z^*\u)\,.\,\mb{n},\varphi>_{\G^-}=0$.
 Thus, we obtain $(z^*\u)\,.\,\mb{n}_{|\G^-}=0$, which is equivalent to $(z\u)\,.\,\mb{n}_{|\G^-}=h$ and the boundary condition (\ref{1b2}) follows.
 \hfill$\diamondsuit$\\[0.3cm]
 \section{Existence of a solution for Problem (\ref{2eqe}) in convex polygon}
 \subsection{Additional regularity in a convex polygon} 
 \H From now on, we assume that $\O$ is a convex polygon. Let $\G_j$, for $1\le j\le N$, denote the sides of $\P\O$, with the convention that $\G_j$ is adjacent to $\G_{j+1}$ and $\G_{N+1}$ coincide with $\G_1$. Also, we denote by $\mb{n}_j$ the corresponding exterior unit normal to $\G_j$, by $\tho_j$ the unit tangent vector along $\G_j$ pointing in the clockwise direction  and by $\mb{x}_j$ the common vertex of $\G_j$ and $\G_{j+1}$.\\
 \H As we shall see later, the regularity $H^1$ of the solution $z$ of the transport equation \ref{1et} requires the regularity $L^{\iy}$ of $\n \u$.
By Sobolev's imbedding theorem, this holds if $\u$ is in
$W^{2,r}(\O)^2$ for some $r>2$. 
 In Proposition 5.3 of [23], they proved
this regularity when $\O$ is a convex polygon, but in the particular case where $\mb{g}\,.\,\mb{n}=0$ on the boundary. The previous lemma establishes this 
regularity without this last assumption. In order to insure the existence of a lifting of $h$, we assume that
\begin{equation}\label{4hh}
h\in W^{1-1/t,t}(\G^-)\H \mr{for}\ t>2.\end{equation}
Then, there exists a lifting $z_h\in W^{1,t}(\O)$ such that
\begin{equation}\label{4dzh}
z_h=h\H\mr{on}\H \G^-\H\mr{with}\ \p z_h\p_{W^{1,t}(\O)}=\p h\p_{W^{1-1/t,t}(\G^-)}.\end{equation}
\begin{lem}\label{4rw2r}
In addition to the hypotheses of Theorem \ref{2te}, we suppose $\O$ is a convex polygon and the boundary data $h$ belongs to $W^{1-1/t,t}(\G^-)$. There exists a real number $r_0>2$,
depending on the inner angles of $\P\O$, such that: if for some real
number $r$ with $2<r<r_0$ and for $1\le j\le N$:
\begin{equation}\label{4crg}
\g_{|\G_j}\in W^{2-1/r,r}(\Gamma_j)^2,\end{equation}
\begin{equation}\label{4ccg}
\g_{|\Gamma_j}(\mb{x}_j)=\g_{|\Gamma_{j+1}}(\mb{x}_j),
\ \F{\P(\g_{|\Gamma_j}\,.\,\mb{n}_{j+1})}{\P \tho_j}(\mb{x}_j)
= \F{\P(\g_{|\Gamma_{j+1}}\,.\,\mb{n}_j)}{\P \tho_{j+1}}(\mb{x}_j)
,\end{equation}then any solution $\u\in W$ of (\ref{1eu})-(\ref{1b2}) belongs
to $W^{2,r}(\O)^2$ (therefore to $W^{1,\iy}(\O)^2$) and
\begin{equation}
\p\u\p_{W^{1,\iy}(\O)}\le C_{\iy}C_r(\sum_{j=1}^N\p\g\p_{W^{2-1/r,r}(\Gamma_j)}
+K(\a,\nu,\f,\g,h)),\label{4bnli}
\end{equation}
where $C_r$ is a constant independent of $\a$ and $\nu$ defined in (\ref{4mu2r}), $C_{\iy}$ is the Sobolev constant defined by (\ref{2dCinf}) and where $K(\a,\nu,\f,\g,h)$ are defined in (\ref{4dK}).
\end{lem}
\tb{Proof.} We cannot prove the regularity $L^{\iy}$ of $\n\u$ by using the generalized Stokes problem (\ref{1esg}), (\ref{1ed}), (\ref{1b1}), because we have only the regularity $L^2$ of $z$. As in [23], we shall use the equality $z=\cs(\u-\a \D\u)$. Since $y=\cs(\D\u)$ belongs to $L^2(\O)$, there exists $\w$ in $H^1(\O)^2$ such that \begin{equation}\label{4rwy} y=\cs\w\H\mr{with}\H \p\w\p_{H^1(\O)}\le C\p y\p_{L^2(\O)}.\end{equation} 
In view of $\cs(\D\u-\w)=0$ and since $\O$ is simply connected, there exists a function $q$ in $L^2(\O)$ such that 
\begin{equation}\label{4stokes}
\D\u-\w=\n q\Longleftrightarrow -\D\u+\n q=-\w,\end{equation}
which implies that the pair $(\u,q)$ is the solution of a Stokes problem with the right-hand side in $H^1(\O)^2$ and its regularity is determined by the angles of $\P\O$ and the regularity of its trace on $\P\O$. In the same way as in [23], the assumptions (\ref{4crg}) and (\ref{4ccg}) (see [1]) imply that there exists a lifting 
$\u_{\g}\in W^{2,r}(\O)^2$ of $\g$ with
$$\d\u_{\g}=0\H\mr{in}\ \O\H\mr{and}\H \p\u_{\g}\p_{W^{2,r}(\O)}\le C\sum_{j=1}^N\p\g\p_{W^{2-1/r,r}(\Gamma_j)}.$$
Therefore the regularity of $\u$ is the same as the regularity of the solution of a homogeneous Stokes problem with the right-hand side in $L^r(\O)^2$ for $r>2$. 
Since $\O$ is a convex polygon and since the right hand-side $\w$ belongs to $L^r(\O)^2$ for all $r\ge 1$, all its inner angles $\o_j$ satisfy $0<\o_j<\pi$ and there exists a real number $r_0>2$ (see Theorem 7.3.3.1 of [25]), depending on the largest inner angle $\o_j$, such that, for some $r<r_0$,  the solution $\u$ belongs to $W^{2,r}(\O)^2$ and, in view of Sobolev imbeddings, to $W^{1,\iy}(\O)^2$, with the existence of a constant $C_{\iy}$ verifying
\begin{equation}\label{2dCinf}
\p\u\p_{W^{1,\iy}(\O)}\le C_{\iy}\p\u\p_{W^{2,r}(\O)}.\end{equation} Hence, owing to (\ref{4rwy}) and (\ref{4stokes}), we have
\begin{equation}\label{4mu2r}
\p\u\p_{W^{2,r}(\O)}\le C_r(\sum_{j=1}^N\p\g\p_{W^{2-1/r,r}(\Gamma_j)}+\p y\p_{L^2(\O)}).\end{equation}
Considering that $\a\, y=\cs\u-z$ and owing to (\ref{2muz2}) with $\varepsilon=\di\F{\nu}{\sqrt{2}\sqrt{S_2^2+1}}$, we obtain
\begin{eqnarray*}|\a|\p y\p_{L^2(\O)}\le \sqrt{2}\F{S_2\sqrt{S_2^2+1}}{\nu}\p\f\p_{L^2(\O)}+2\p z\p_{L^2(\O)}\\+\sqrt{2}(1+\sqrt{S_2^2+1}\,) L(\F{\nu}{\sqrt{2}\sqrt{S_2^2+1}\p\g\p_{H^{1/2}(\P\O)}})\p\g\p_{H^{1/2}(\P\O)}.\end{eqnarray*}
Substituting the estimate (\ref{2mz}) into this last inequality, considering (\ref{2dCinf}) and using the bound (\ref{4mu2r}) yield (\ref{4bnli}) with 
\begin{eqnarray}
K(\a,\nu,\f,\g,h)=\tilde{C}(\a,\nu,\f,\g,h)+\F{2}{|\a|}C(\a,\nu,\f,\g,h),\label{4dK}\end{eqnarray} where
\begin{eqnarray*}\tilde{C}(\a,\nu,\f,\g,h)=\F{\sqrt{2}}{\nu|\a|}\sqrt{S_2^2(S_2^2+1)+8\a^2}\p\f\p_{H(\cs,\O)}\\ +\F{8}{|\a|}\p z_h\p_{L^2(\O)} +\F{\sqrt{2}(1+\sqrt{S_2^2+1}\,)}{|\a|} L(\F{\nu}{\sqrt{2}\sqrt{S_2^2+1}\p\g\p_{H^{1/2}(\P\O)}})\p\g\p_{H^{1/2}(\P\O)}\end{eqnarray*} and 
$C(\a,\nu,\f,\g,h)$ is defined by (\ref{3dC}), which ends the proof of the lemma.\hfill$\diamondsuit$\\[0.2cm]
\subsection{Existence of a solution for Problem (\ref{2eqe})}
\H In this section, the domain $\O$ is a convex polygon and with adequate asumptions on $\g$ and $h$, we shall see that the previous solutions $z$ will belong to
$H^1(\O)$. First, we define some notations that we will need to specify certain assumptions of the next theorem. 
Let $\x$ belong to $\gamma_+\cap\gamma_-$ or to $\stackrel{\circ}{\gamma_-}$, where $\gamma_+$ and $\gamma_-$ are two straight segments such that $\gamma_+\subset \ov{\G^{+,0}}$ and $\gamma_-=[\x,\mb{x_-}]\subset \ov{\G^-}$. We denote by
 \begin{equation}\label{4dnm}
 \mb{n_-(\x)}\ \mr{the\ unit\ exterior\ normal\ vector\ to}\ \gamma_-,\end{equation} by 
 \begin{equation}\label{4dtom}
 \tho_-(\x)\ \mr{the\ unit\ tangent\ vector}\ \F{1}{\p \x\mb{x_-}\p}\x\mb{x_-}\end{equation} and by 
 $E$ the set defined by
 \begin{equation}\label{4dE}
E=\{\x\in \ov{\Gamma^-}\cap\ov{\Gamma^{+,0}},\ \g(\x)\,.\,\mb{n}_-(\x)=0\}.\end{equation}
Note that, in view of the assumption (\ref{1dg0g1}), the set $E$ is finite. In addition, we make the assumption that the data $\g$ is such that
\begin{equation}\label{4hs}
\{\x\in\ov{\Gamma^-},\ \g(\x)\,.\,\mb{n}(\x)=0\}\subset E,\end{equation}
where $\mb{n}(\x)$ is the unit exterior normal vector to the boundary at the point $\x$, if $\x$ is not a vertex, and a unit exterior normal vector to $\ov{\Gamma^-}$ at the point $\x$, if $\x$ is a vertex. Note that (\ref{4hs}) means that $\g\,.\,\mb{n}$ does not vanish in a point located in the interior of $\ov{\G^-}$. The following theorem gives assumptions implying existence for Problem (\ref{2eqe}).
\begin{theo}\label{4te} Let $\O$ be a convex polygon, let $\f$ be in $L^2(\O)^2$ such that $\cs\f\in H^1(\O)$, let the real numbers $\nu$, $\a$ and $t$ such that $\nu>0$ and 
$t>2$, let $h$ belong to $W^{1-1/t,t}(\G^-)$ and let $\g$ be such that $\g_{|\G_j}\in W^{2-1/r,r}(\G_j)^2$ for some real $r$ with $2<r<r_0$ and for $1\le j\le N$, satisfying the second part of (\ref{1b1}) and (\ref{4hs}) and verifying (\ref{4crg}), (\ref{4ccg}) and 
\begin{equation}\label{4hzH1}
\forall \x\in E,\ \F{\P \g}{\P \tho_-}(\x)\,.\,\mb{n_-(\x)}\not=0\ \mr{and}\ \g(\x)\,.\,\tho_-(\x)<0,\end{equation}
where $E$, $\mb{n_-}$ and $\tho_-$ are defined in (\ref{4dE}), (\ref{4dnm}) and (\ref{4dtom}).
We suppose, in addition, that the data $\a$, $\nu$, $\f$, $\g$ and $h$ are small such that
\begin{equation}
 C_{\iy}C_r(\sum_{j=1}^N\p\g\p_{W^{2-1/r,r}(\Gamma_j)}
+K(\a,\nu,\f,\g,h))\le \F{\nu}{2|\a|}\label{4czH1},\end{equation}
where $C_r$ is a constant independent of $\a$ and $\nu$ defined by (\ref{4mu2r}), $C_{\iy}$ is the Sobolev constant defined by (\ref{2dCinf}) and where the function $K$  is defined in (\ref{4dK}). Then, there exists at least one solution $(\u,p,z)$ for Problem (\ref{2eqe}).\end{theo}
\tb{Proof.} First, with a new definition of $z_h$ (defined by (\ref{4dzh}) instead of (\ref{2nzh})), we define a sequence $(z^*_n)$ in the same way as in (\ref{3dzn}), (\ref{3eqe}) and (\ref{3eqt}) and by the same convergenge as in the proof of Theorem \ref{2te}, we obtain a solution $(\u,p,z)\in (H^1(\O))^2\t L^2_0(\O)\t L^2(\O)$, where $z=z^*+z_h$, that verifies
\begin{eqnarray}
 -\nu\D\u+\z\t\u+\n p=\f\H\mr{in}\ \O,\nonumber\\
  \d\,\u=0\hspace{1cm} \mathrm{in}\hspace{0.3cm}\O,\nonumber\\
\u=\g\hspace{0.8cm} \mr{on}\ \P\O,\label{4eqe}\\
  \nu z +\a \u\,.\,\n z=\nu\,\cs\u+\a\,\cs\f\H\mr{in}\ \O,\nonumber\\
 (z\u)\,.\,\mb{n}=(h\u)\,.\,\mb{n}\H\mr{on}\ \G^-.\nonumber
 \end{eqnarray} 
 This solution $(\u,p,z)$, owing to Lemma \ref{1leqf}, is solution of Problem $(\mr{P_I})$, that is to say of the system (\ref{1eu})-(\ref{1b2}). Applying Lemma 
 \ref{4rw2r} yields that the velocity $\u$ belongs to $W^{2,r}(\O)^2$, therefore to $W^{1,\iy}(\O)^2$, and satifies the estimate (\ref{4bnli}). In view of the assumption 
 (\ref{4czH1}), we derive that the velocity $\u$ verifies the bound
 \begin{equation}\label{4ubwi}
 |\a|\,\p\u\p_{W^{1,\iy}(\O)}\le \F{\nu}{2}.
 \end{equation}
 Next, with this velocity $\u$, we associate the following problem : find $\tilde{z}\in H^1(\O)$ solution of the transport problem
 \begin{equation}\label{4eqt}\left\lbrace\begin{array}{ll}
 \nu\tilde{z}  +\a \u\,.\,\n \tilde{z}=
 \nu\,\cs\u+\a\,\cs\f\H&\mr{in}\ \O,\\
 \tilde{z}=h&\mr{on}\ \G^-.\end{array}\right.\hspace*{1.2cm}  
 \end{equation}
If we show, first, that this problem has a unique solution $\tilde{z}$ and, second, that $\tilde{z}=z$, where $(u,p,z)$ is the previous solution of (\ref{4eqe}), then we will have proven the existence of a solution $(\u,p,z)$ for Problem (\ref{2eqe}). Using $z_h$ defined by (\ref{4dzh}), let us split the problem (\ref{4eqt}) into two transport problems, namely: find $(\tilde{z_1},\tilde{z_2})\in H^1(\O)^2$ such that
\begin{equation}\label{4eqt1}\left\lbrace\begin{array}{ll}
 \nu\tilde{z_1}  +\a \u\,.\,\n \tilde{z_1}=
 \nu\,\cs\u+\a\,\cs\f-\nu z_h\H&\mr{in}\ \O,\\
 \tilde{z_1}=0&\mr{on}\ \G^-\end{array}\right.
 \end{equation}
 and   
 \begin{equation}\label{4eqt2}\left\lbrace\begin{array}{ll}
 \nu\tilde{z_2}  +\a \u\,.\,\n \tilde{z_2}=
 \nu z_h\H&\mr{in}\ \O,\\
 \tilde{z_2}=h&\mr{on}\ \G^-.\end{array}\right.\hspace*{3.55cm} 
 \end{equation}
 The first problem is a transport problem from the type
 $$z +\mathcal{W} \u\,.\,\n z=l,\ \mr{with}\ z_{|\G^-}=0$$
 with $\mathcal{W}=\di\F{\a}{\nu}$ and $l=\,\cs\u+\F{\a}{\nu}\,\cs\f- z_h$ and this problem was solved by Theorem 3.1 of [7], the assumptions of which are verified owing to (\ref{4hzH1}) and since 
 (\ref{4ubwi}) implies to $\p\n\u\p_{L^{\iy}(\O)}\le \di\F{1}{2|\mathcal{W}|}$. Thus, we obtain a solution $\tilde{z_1}\in H^1(\O)$ of (\ref{4eqt1}).\\
 \H The second problem is a little different because of the non-homogeneous boundary condition on $\G^-$. In order to solve this problem, we shall use an analogous method as in Theorem 2.1 of [7]. We define  a sequence ($\mb{F}_n$) of function $\mb{F}_n\in X_{\u}(\G^-)^2$, $n\in \N$, by recurrence.
  We set $\mb{F}_0=\mb{0}$ and assume that the function $\mb{F}_n\in X_{\u}(\G^-)^2$ is given for $n\in\N$. Then, applying Theorem 3.3 of [5], we define each component $F_{n+1,1}$ and $F_{n+1,2}$ of $\mb{F}_{n+1}$ as the unique solution of the transport problem from the type $z+\mathcal{W}\,\u\,.\,\n z=l$ with $(z\u)\,.\,\mb{n}=0$, so that we define $\mb{F}_{n+1}\in X_{\u}(\G^-)^2$ as the unique solution of the transport problem
\begin{equation}\left\lbrace\begin{array}{ll}
\mb{F}_{n+1}+\di\F{\a}{\nu}\,\u\,.\,\n \mb{F}_{n+1}=\n z_h-\di\F{\a}{\nu}\,\n\u\,.\,\mb{F}_{n}\H &\mr{in}\H \O\\[0.2cm]
(\mb{F}_{n+1}\,\u)\,.\,\mb{n}=0 &\mr{on}\H \G^-.\end{array}\right.\label{4eqt3}
\end{equation}
Since $\mb{F}_{n+1}$ belongs to $X_{\u}(\G^-)^2$, the basic result of Proposition (\ref{2pgreenfg-}) implies
$$\int_{\O}(\F{\a}{\nu}\,\u\,.\,\n F_{n+1,i})F_{n+1,i}\,d\mb{x}\ge 0,\ \mr{for}\ i=1,2.$$
 Then, taking the scalar product of both sides of the first equation of (\ref{4eqt3}) with $\mb{F}_{n+1}$ yields
$$\p\mb{F}_{n+1}\p_{L^2(\O)}^2\le (\n z_h,\mb{F}_{n+1})-\F{\a}{\nu}\,(\n\u\,.\,\mb{F}_{n},\mb{F}_{n+1}).$$
In view of the bound (\ref{4ubwi}), we obtain
$$\p\mb{F}_{n+1}\p_{L^2(\O)}\le \p\n z_h\p_{L^2(\O)}+\F{1}{2}\p\mb{F}_{n}\p_{L^2(\O)},$$
which implies, by a recurrence argument, that $\mb{F}_{n}$ is uniformly bounded in $L^2(\O)$ and 
$\forall n\in \N$,
\begin{equation}\label{4bFn}
\p\mb{F}_{n}\p_{L^2(\O)}\le 2\p\n 
z_h\p_{L^2(\O)}.\end{equation}
Owing to (\ref{4bFn}), $\u\,.\,\n \mb{F}_{n+1}$ is also uniformly bounded in $L^2(\O)$. Therefore we can pass to the limit in the first equation of 
(\ref{4eqt3}) and there exists a function $\mb{F}\in L^2(\O)^2$ such that
\begin{equation}\label{4eqtd}\mb{F}+\di\F{\a}{\nu}\,(\u\,.\,\n\mb{F}+\n\u\,.\,\mb{F})=\n z_h.\end{equation}
Let us set $\tilde{z_2}=z_h-\di\F{\a}{\nu}\u\,.\mb{F}$. From the previous equation, we derive $\mb{F}=\n \tilde{z_2}$ and we obtain $\tilde{z_2}=z_h-\F{\a}{\nu}\u\,.\,\n \tilde{z_2}$, which gives that $\tilde{z_2}\in H^1(\O)$ is solution the first equation of (\ref{4eqt2}).\\
In the same way as in the proof of Theorem 2.1 of [7], we prove 
$$(\u\,.\,\n\tilde{z_2})_{|\G^-}=0.$$
Hence, considering that $\nu\, \tilde{z_2}+\a\u\,.\,\n\tilde{z_2}=\nu\, z_h$, we obtain 
$$\tilde{z_2}_{|\G^-}=(z_h)_{|\G^-}=h,$$
which implies that $\tilde{z_2}$ is solution of (\ref{4eqt2}) and verifies
\begin{equation}\label{4mnz2}
\p\n z_2\p_{L^2(\O)}\le  2\p\n 
z_h\p_{L^2(\O)}.\end{equation} Finally,  $$\tilde{z}=\tilde{z_1}+\tilde{z_2}$$
is a solution of (\ref{4eqt}). Let us show that $\tilde{z}$ also verifies the boundary condition on $\G^-$ of (\ref{4eqe}). $ \forall\vfy\in W^{1-1/r,r}_{00}(\G^-),$
$$<(\tilde{z}\u)\,.\,\mb{n},\vfy>_{\G^-}=\int_{\G^-}(\tilde{z}\vfy \u)\,.\,\mb{n}\,ds= \int_{\G^-}(h\vfy \u)\,.\,\mb{n}\,ds=<(h\u)\,.\,\mb{n},\vfy>_{\G^-},$$
where $W^{1-1/r,r}_{00}(\G^-)$ is defined by (\ref{1dW00r}) with $r>2$. Hence, we derive $$(\tilde{z}\u)\,.\,\mb{n}_{|\G^-}=(h\u)\,.\,\mb{n}_{|\G^-},$$ which implies that 
$z$ and $\tilde{z}$ verify the same transport equation with the same boundary condition. Then, $Z=z-\tilde{z}$ is solution of the transport problem
\begin{equation}\label{4eqt4}\left\lbrace\begin{array}{ll}
 \nu\, Z  +\a \u\,.\,\n Z=0
 \H&\mr{in}\ \O,\\
 (Z\u)\,.\,\mb{n}=0&\mr{on}\ \G^-.\end{array}\right.
 \end{equation}
 Owing to the Proposition (\ref{2pgreenfg-}), which implies $( \a \u\,.\,\n Z,Z)\ge 0$, we have $Z=0\Longleftrightarrow z=\tilde{z}$. Finally, $z\in H^1(\O)$ and verify $z_{|\G^-}=h$, so $(\u,p,z)$ is solution of Problem (\ref{2eqe}).\hfill$\diamondsuit$\\[0.2cm]
\section{Uniqueness}
\H In order to establish uniqueness of the solution of the fully non-homogeneous problem of grade 2 fluids, we need the regularity $H^1$ of the solution $z$ of the transport equation (\ref{1et}), which is the framework of the Problem \ref{2eqe}. Moreover, if we want conditions of uniqueness independent from $z$, we need a bound of this solution in $H^1(\O)$. However, to get this bound,
we are led to make more restrictive assumptions as in the previous theorem.
\subsection{$H^1$ Bound for the transport equation}
\H We deal with a transport equation from the type
\begin{equation}\label{5eqt}\left\lbrace\begin{array}{ll}
 z  +\mathcal{W} \u\,.\,\n z=
 l\H&\mr{in}\ \O,\\
 z=h&\mr{on}\ \G^-.\end{array}\right.\hspace*{3.55cm} 
 \end{equation}
 \H If we make the restrictive assumption that $\u\,.\,\mb{n}$ does not vanish at the boundary of $\G^-$, as in Theorem 2.2 of [7], then we define the real number $\beta>0$ by
 \begin{equation}\label{5dbeta}
 \beta=\max (\F{1}{\u(\x)\,.\,\mb{n}_-(\x)},\ \x\in\ov{\G^-}).\end{equation}
 Note that the assumptions of Theorem 3.1 in [7] are less restrictive and under these assumptions, we can prove the existence of a solution $z$ in $H^1$ of the transport equation  (\ref{5eqt}) in the particular case where $h=0$. Yet, in this frame, we are not able to have a bound in $H^1$ of this solution. In fact, we can bound the solution of (\ref{5eqt}) in $H^1$ only in the case where $l_{|\G^-}=h$. In the same way, we can bound the solution $\tilde{z_2}$ of 
 (\ref{4eqt2}) in $H^1$, but we cannot bound the solution $\tilde{z_1}$ of (\ref{4eqt1}). Unhappily, if we want conditions of uniqueness only depending on the data, the uniqueness of the solution of the fully non-homogeneous grade 2 problem seems to require a bound of the solution $z$ in $H^1$, which leads to impose the restrictive assumption that $\u\,.\,\mb{n}$ does not vanish on $\ov{\Gamma^-}$.
  \begin{theo}\label{5thbH1} Let $\O$ be a bounded polygon, $\G^-$ be defined by (\ref{1dg-}), verifying 
(\ref{1dg0g1}). For all $\u$ in 
 $W^{1,\iy}(\O)^2$, verifying $\d\u=0$,  such that
 \begin{equation}\label{5ctheu}\p\n\u\p_{L^{\iy}(\O)}\le \di\F{1}{2|\W|}\end{equation} and such that
 \begin{equation}\label{5cun}
 \forall \x\in\ov{\Gamma^-},\ \u(\x)\,.\,\mb{n}_-(\x)\not=0,\end{equation} all $l$ in $H^1(\O)$, all $h$ in $H^{1/2}(\G^-)$
  and all real number $\W$ in $\R^*$, the transport problem (\ref{5eqt}) has a unique solution $z$ in $H^1(\O)$ and this solution satisfies the following estimate:
  \begin{equation}\label{5mzH1}
  \p\n z\p_{L^2(\O)}\le (2+\beta C_0(\F{3}{|\mathcal{W}|}+2\p\u\p_{W^{1,\iy}(\O)}))(\p l\p_{H^1(\O)}+\p h\p_{H^{1/2}(\G^-)}),
  \end{equation}
  where $C_0$ is defined by (\ref{5dC0}) and $\beta$ by (\ref{5dbeta}).
\end{theo}
\tb{Proof.} In view of (\ref{5cun}), the regularity of $l$, $\u$ and $h$ imply that $\di\F{l_{|\G^-}-h}{\W (\u\,.\mb{n})_{|\G^-}}$ belongs to $H^{1/2}(\G^-)$. So, 
there exists a lifting $z_0\in H^2(\O)$ and a real constant $C_0$ such that
\begin{equation}\label{5dy0}\left\lbrace\begin{array}{ll}
 \di\F{\P z_0}{\P n}_{|\G^-}=\di\F{l_{|\G^-}-h}{\W\,\u\,.\,\mb{n}}\\[0.2cm]
 z_0\hspace*{0.01cm}_{|\G^-}=0
 \end{array}\right. 
 \end{equation} and 
 \begin{equation}\label{5dC0}
 \p z_0\p_{H^2(\O)}\le C_0 \p\di\F{l_{|\G^-}-h}{\W\,\u\,.\,\mb{n}}\p_{H^{1/2}(\G^-)}.\end{equation}
\H Next, we define the following problem : find $z^*$ in $H^1(\O)$  solution of the transport problem
\begin{equation}\label{5eqtz*}\left\lbrace\begin{array}{ll}
 z^*  +\mathcal{W} \u\,.\,\n z^*=l^*
 \H&\mr{in}\ \O,\\
 z^*=h&\mr{on}\ \G^-,\end{array}\right.\hspace*{3.55cm} 
 \end{equation}
 where 
 \begin{equation}\label{5dl*}l^*=l-z_0-\W\,\u\,.\,\n z_0.\end{equation}
 Note that, by construction of $z_0$, we have $l^*_{|\G^-}=h$. Hence, we derive that 
 $$(\u\,.\,\n z^*)_{|\G^-}=0.$$  
  In the same way as in (\ref{4eqt3}), we define a sequence $(\mb{F}^*_n)$, $\mb{F}^*_n\in X_{\u}(\G^-)^2$, by 
 \begin{equation}\left\lbrace\begin{array}{ll}
\mb{F}_{n+1}^*+\W\,\u\,.\,\n \mb{F}_{n+1}^*=\n l^*-\W\,\n\u\,.\,\mb{F}_{n}^*\H &\mr{in}\H \O\\[0.2cm]
(\mb{F}_{n+1}^*\,\u)\,.\,\mb{n}=0 &\mr{on}\H \G^-.\end{array}\right.\label{5eqt3*}
\end{equation}
By the same method as previously we prove that
$$\p\mb{F}_n^*\p_{L^2(\O)}\le 2\p \n l^*\p_{L^2(\O)}$$ and the sequence $(\mb{F}^*_n)$ converge to $\mb{F}^*\in X_{\u}(\O)^2$, which verifies 
\begin{equation}\label{5mF*}
\mb{F}^*+\W\u\,.\,\n\mb{F}^*  =\n l^*-\W\,\n\u\,.\,\mb{F}^*\H \mr{and}\H \p\mb{F}^* \p_{L^2(\O)}\le 2 \p \n l^*\p_{L^2(\O)}.\end{equation}
 Then, setting 
 \begin{equation}\label{5dz*}
  z^*=l^*-\W \u\,.\,\mb{F}^*,
 \end{equation} we obtain that $\mb{F}^*=\n z^*$ and, since we can prove that $(\u\,.\,\mb{F}^*)_{|\G^-}=0$ as in the proof of Theorem 2.1 of [7], $z^*$ is the solution  of Problem (\ref{5eqtz*}).
 Finally
 $$
 z=z^*+z_0$$ is the unique solution of (\ref{5eqt}). Hence, from (\ref{5mF*}) and (\ref{5dl*}), in view of $\mb{F}^*=\n z^*$, we derive
 \begin{eqnarray}
 \p\n z\p_{L^2(\O)}\le 3\p \n z_0\p_{L^2(\O)}+2\p\n l\p_{L^2(\O)}+2|\W|\,\p\u\p_{W^{1,\iy}(\O)}\p z_0\p_{H^2(\O)}\nonumber\\
 \le 2\p\n l\p_{L^2(\O)}+(3+2|\W|\,\p\u\p_{W^{1,\iy}(\O)})\p z_0\p_{H^2(\O)}.\hspace*{1.52cm}\label{5mz*i}\end{eqnarray}
 Next, owing to (\ref{5dC0}) and (\ref{5dbeta}), we obtain
 $$\p z_0\p_{H^2(\O)}\le \F{C_0\beta}{|\W|}(\p l\p_{H^1(\O)}+\p h\p_{H^{1/2}(\G^-)}).$$
 Substituting this bound into (\ref{5mz*i}) yields (\ref{5mzH1}).\hfill$\diamondsuit$\\[0.2cm]
 \subsection{Uniqueness}
 \H The transport equation $\nu\, z+\a\,\u\,.\,\n z=\nu\,\c\u+\a\,\c\f$ is of the type\linebreak $z+\W\u\,.\,\n z=l$ with $\W=\di\F{\a}{\nu}$ and $l=\c\u+\di\F{\a}{\nu}\c\f$. In the 
 estimate (\ref{5mzH1}), the norm $H^1$ of $l$ occurs, so we need to bound the norm $H^2$ of $\u$. This bound is given in the following lemma the proof of which is analogous to that of Lemma \ref{4rw2r} and where the constant $C_2$ is defined by
 \begin{equation}\label{4muH2}
\p\u\p_{H^2(\O)}\le C_2(\sum_{j=1}^N\p\g\p_{H^{3/2}(\Gamma_j)}+\p y\p_{L^2(\O)}).\end{equation}
  \begin{lem}\label{4rw2}
In addition to the hypotheses of Theorem \ref{2te}, we suppose $\O$ is a convex polygon, the boundary data $h$ belongs to $W^{1-1/t,t}(\G^-)$ and for $1\le j\le N$, the boundary data $\g$ verifies
\begin{equation}\label{4c2g}
\g_{|\G_j}\in H^{3/2}(\Gamma_j)^2\end{equation} and (\ref{4ccg}), then any solution $\u\in W$ of (\ref{1eu})-(\ref{1b2}) belongs
to $H^2(\O)^2$  and
\begin{equation}
\p\u\p_{H^2(\O)}\le C_2(\sum_{j=1}^N\p\g\p_{H^{3/2}(\Gamma_j)}
+K(\a,\nu,\f,\g,h)),\label{5buH2}
\end{equation}
where $C_2$ is a constant independent of $\a$ and $\nu$ defined by (\ref{4muH2}) and where the function $K$ is defined in (\ref{4dK}).
\end{lem}
\H Under the assumptions of Theorem \ref{4te}, Problem (\ref{2eqe}) has at least one solution $(\u,p,z)$ and, owing to Theorem \ref{5thbH1} and the previous lemma, we can bound the function $\n z$ in $L^2(\O)$. Note that, since $\u=\g$ on $\G^-$, we can express the constant $\beta$, defined by (\ref{5dbeta}), by using $\g$ as
\begin{equation}\label{5dbetag}
\beta=\max (\F{1}{\g(\x)\,.\,\mb{n}_-(\x)},\ \x\in\ov{\G^-}).\end{equation}
\begin{lem}\label{5lbzH1bis} Under the assumptions of Theorem \ref{4te}, let $(\u,p,z)$ a solution of Problem (\ref{2eqe}). Then, we have the following bound 
\begin{equation}\label{5bzH1bis}
\p\n z\p_{L^2(\O)}\le L_3(\a,\nu,\f,\g,h),\end{equation}
with \begin{equation}\label{5dL3}
L_3(\a,\nu,\f,\g,h)=\tilde{L}_3(\a,\nu,\f,\g,h)\hat{L}_3(\a,\nu,\f,\g,h),\end{equation}
where the functions $\tilde{L}_3$ and $\hat{L}_3$ are defined by 
$$
\tilde{L}_3(\a,\nu,\f,\g,h)=2+\F{3 C_0\beta\nu}{|\a|}+2\beta C_0 C_{\iy}C_r(\sum_{j=1}^N\p\g\p_{W^{2-1/r,r}(\Gamma_j)}
+K(\a,\nu,\f,\g,h)),$$
$$\hat{L}_3(\a,\nu,\f,\g,h)=\p h\p_{H^{1/2}(\G^-)}+\F{|\a|}{\nu}\p\c\f\p_{H^1(\O)}+\sqrt{2} C_2(\sum_{j=1}^N\p\g\p_{H^{3/2}(\Gamma_j)}
+K(\a,\nu,\f,\g,h)),$$
where $C_0$ is defined by (\ref{5dC0}), $\beta$ by (\ref{5dbetag}), $C_{\iy}$ by (\ref{2dCinf}), $C_r$ by (\ref{4mu2r}), $C_2$ by (\ref{4muH2}) and where the function $K$ are defined in (\ref{4dK}). 
\end{lem}
\tb{Proof.} Applying Theorem \ref{5thbH1} with $\W=\di\F{\a}{\nu}$ and $l=\c\u+\di\F{\a}{\nu}\c\f$ we obtain
$$
\p\n z\p_{L^2(\O)}\le(2+ 3\F{C_0\beta\nu}{|\a|}+2\beta C_0\p\u\p_{W^{1,\iy}(\O)})(\sqrt{2}\p\u\p_{H^2(\O)}+\F{|\a|}{\nu}\p\c\f\p_{H^1(\O)}+\p h\p_{H^{1/2}(\G^-)}).$$
Substituting the bounds (\ref{4bnli}) and (\ref{5buH2}) in this last inequality yields (\ref{5bzH1bis}).\hfill$\diamondsuit$\\[0.2cm]
\H Now, we give two results of uniqueness for Problem (\ref{2eqe}). Note that, in Theorem \ref{5thupb2}, we do not suppose that $E$ is an empty set, but the condition (\ref{5cuz}) depends on the norm $\p\n z\p_{L^2(\O)}$ where $(\u,p,z)$ is any solution of Problem (\ref{2eqe}), while, in  Theorem \ref{5thupb2b}, we suppose that $E=\emptyset$, that is to say the condition 
(\ref{5cun}) with $\u=\g$ on $\G^-$, which is a more restrictive condition that in Theorem \ref{5thupb2}, the hypotheses of which follows those of Theorem \ref{4te}, which insures the existence of solutions. Indeed, without 
this restrictive condition, we are not able to bound the norm $\p\n z\p_{L^2(\O)}$ in a function of the data. In conclusion, contrary to 
Theorem \ref{5thupb2}, the conditions of uniqueness of Theorem \ref{5thupb2b} only depend on the data, but its conditions of uniqueness are more restrictive.
\begin{theo}\label{5thupb2} Let $\O$ be a convex polygon, let $\f$ be in $L^2(\O)^2$ such that $\cs\f\in H^1(\O)$, let the real numbers $\nu$, $\a$ and $t$ such that $\nu>0$ and 
$t>2$, let $h$ belong to $W^{1-1/t,t}(\G^-)$, let $\g$ be such that $\g_{|\G_j}\in W^{2-1/r,r}(\G_j)^2$ for some real $r$ with $2<r<r_0$ and for $1\le j\le N$, satisfying the second part of (\ref{1b1}) and (\ref{4hs}) and verifying (\ref{4crg}), (\ref{4ccg}) and (\ref{4hzH1}.
We suppose, in addition, that the data $\a$, $\nu$, $\f$, $\g$ and $h$ are small such that they verify (\ref{4czH1}) and such that any solution $(\u,p,z)$ of Problem (\ref{2eqe}) satisfies
\begin{equation}
\F{S^*_4}{\nu}L_2(\a,\nu,\f,\g,h)(\sqrt{2} S_4+\F{|\a|}{\nu}C_{3/2}S_{\iy}\bar{S}(1+\F{S_4^2}{\nu} L_1(\a,\nu,\f,\g,h))\p\n z\p_{L^2(\O)}< 1
\label{5cuz},\end{equation}
where $C_r$ and $C_{3/2}$ are constants independent of $\a$ and $\nu$ defined respectively by (\ref{4mu2r}) and (\ref{5dC3/2}), $C_{\iy}$, $S_4$, $S^*_4$, $S_{\iy}$ and $\bar{S}$ are Sobolev constants defined respectively by (\ref{2dCinf}), (\ref{1dSp}), (\ref{1dSp*}), (\ref{5dSinf}) and (\ref{5dbarS}) and where the functions $K$, $L_1$ and $L_2$ are defined respectively in (\ref{4dK}), (\ref{5dL1}) and (\ref{5dL2}). Then Problem (\ref{2eqe}) has a unique solution
.\end{theo}
\begin{theo}\label{5thupb2b}
Let $\O$ be a convex polygon, let $\f$ be in $L^2(\O)^2$ such that $\cs\f\in H^1(\O)$, let the real numbers $\nu$, $\a$ and $t$ such that $\nu>0$ and 
$t>2$, let $h$ belong to $W^{1-1/t,t}(\G^-)$ and let $\g$ be such that $\g_{|\G_j}\in W^{2-1/r,r}(\G_j)^2$ for some real $r$ with $2<r<r_0$ and for $1\le j\le N$, satisfying the second part of (\ref{1b1}) and (\ref{4hs}) and verifying (\ref{4crg}), (\ref{4ccg}) and 
\begin{equation}\label{5cgn}
 \forall \x\in\ov{\Gamma^-},\ \g(\x)\,.\,\mb{n}_-(\x)\not=0,\end{equation}
where $\mb{n_-}$ is defined in (\ref{4dnm}).
We suppose, in addition, that the data $\a$, $\nu$, $\f$, $\g$ and $h$ are small such that they verify (\ref{4czH1})
and 
\begin{equation}
\F{S^*_4}{\nu}L_2(\a,\nu,\f,\g,h)(\sqrt{2} S_4+\F{|\a|}{\nu}C_{3/2}S_{\iy}\bar{S}(1+\F{S_4^2}{\nu} L_1(\a,\nu,\f,\g,h))L_3(\a,\nu,\f,\g,h)< 1
\label{5cuzb},\end{equation}
where $C_r$ and $C_{3/2}$ are constants independent of $\a$ and $\nu$ defined respectively by (\ref{4mu2r}) and (\ref{5dC3/2}), $C_{\iy}$, $S_4$, $S^*_4$, $S_{\iy}$ and $\bar{S}$ are Sobolev constants defined respectively by (\ref{2dCinf}), (\ref{1dSp}), (\ref{1dSp*}), (\ref{5dSinf}) and (\ref{5dbarS}) and where the functions $K$, $L_1$, $L_2$ and $L_3$  are defined respectively in (\ref{4dK}), (\ref{5dL1}), (\ref{5dL2}) and (\ref{5dL3}). 
Then, there exists a unique solution $(\u,p,z)$ for Problem (\ref{2eqe}).\end{theo}
\tb{Proof of Theorems \ref{5thupb2} and \ref{5thupb2b}.}\\
\H Let us, first, prove the following lemma.
\begin{lem} Let $(\u_1,p_1,z_1)$ and $(\u_2,p_2,z_2)$ two solutions of Problem (\ref{2eqe}). Then $\u_2-\u_1$ verifies the following estimates :
\begin{equation}\label{5bH1u2u1}
|\u_2-\u_1|_{H^1(\O)}\le\F{S_4}{\nu}\p \u_1\p_{L^4(\O)}\p z_2-z_1\p_{L^2(\O)}\end{equation}
and
\begin{equation}\label{5bliu2u1}
\p \u_2-\u_1\p_{L^{\iy}(\O)}\le\F{C_{3/2} S_{\iy} \bar{S}}{\nu}\p \u_1\p_{L^4(\O)}(1+\F{S_4^2}{\nu}\p z_2\p_{L^2(\O)})\p z_2-z_1\p_{L^2(\O)},
\end{equation}
where $C_{3/2}$ is a constant independent of $\a$ and $\nu$ defined  by (\ref{5dC3/2}), $S_4$, $S_{\iy}$ and $\bar{S}$ are Sobolev constants defined respectively by (\ref{1dSp}), (\ref{5dSinf}) and (\ref{5dbarS}).
\end{lem}
\tb{Proof.} Since $(\u_1,p_1)$ and $(\u_2,p_2)$ are solutions of the generalized Stokes problem (\ref{1esg}), (\ref{1ed}),
 (\ref{1b1}), $\u_2-\u_1$ verifies the following equation
 \begin{equation}\label{5eu2-u1} 
 -\nu\D(\u_2-\u_1)+\n (p_2-p_1)=-\z_2\t (\u_2-\u_1)-(\z_2-\z_1)\t \u_1.\end{equation}
 Then, Green's formula with $\u_2-\u_1$ in $H^1_0(\O)$ yields
 $$\nu|\u_2-\u_1|^2_{H^1(\O)}=-((\z_2-\z_1)\t\u_1,\u_2-\u_1),$$
 which implies
 $$\nu|\u_2-\u_1|^2_{H^1(\O)}\le \p \u_1\p_{L^4(\O)}\p \u_2-\u_1\p_{L^4(\O)}\p z_2-z_1\p_{L^2(\O)}.$$
 Hence, in view of $\p \u_2-\u_1\p_{L^4(\O)}\le S_4|\u_2-\u_1|_{H^1(\O)}$, we obtain (\ref{5bH1u2u1}).\\
 \H Next, considering (\ref{5eu2-u1}),  we derive that $\u_2-\u_1$ is the solution of the homogeneous Stokes's problem with
    $-\z_2\t (\u_2-\u_1)-(\z_2-\z_1)\t \u_1$ as the right-hand side. Note that the right-hand side $-\z_2\t (\u_2-\u_1)-(\z_2-\z_1)\t \u_1$ belongs to $L^{4/3}(\O)\subset H^{-1/2}(\O)$. Applying Theorem 5.5 page 83 of [14], we derive that, for $\theta<1$, if the right-hand side
   belongs to $H^{-1+\theta}$, the homogeneous solution belongs to
   $H^{1+\theta}(\O)^2$, under the
   condition : $$\theta<\min_{(i,j)\in\mathcal{N}_{\O}}\xi_1(\omega_j^i).$$
   But, considering that
   $\forall (i,j)\in \mathcal{N}_{\O}$,
   $\xi_1(\omega_j^i)>\F{1}{2}$, we can apply this result with $\theta=\F{1}{2}$ and we obtain that $\u_2-\u_1$ belongs to $H^{3/2}(\O)^2$  and there exists a positive constant $C_{3/2}$ independent of $\a$ and $\nu$ such that
   \begin{equation}\label{5dC3/2}
   \nu\p\u_2-\u_1\p_{H^{3/2}(\O)}\le C_{3/2}\p-\z_2\t (\u_2-\u_1)-(\z_2-\z_1)\t \u_1\p_{H^{-1/2}(\O)}.\end{equation}
   Then, Sobolev imbeddings, $H^{3/2}(\O)\subset L^{\iy}(\O)$ and $L^{4/3}(\O)\subset H^{-1/2}(\O)$, yield that there exists positive constants $S_{\iy}$ and $\bar{S}$ such that
   \begin{equation}\label{5dSinf}
   \forall v\in H^{3/2}(\O), \p v\p_{L^{\iy}(\O)}\le S_{\iy}\p v\p_{H^{3/2}(\O)},\end{equation}
   \begin{equation}\label{5dbarS}
   \forall v\in H^{-1/2}(\O), \p v\p_{H^{-1/2}(\O)}\le \bar{S}\p v\p_{L^{4/3}(\O)}.\end{equation}
   Hence, we derive
   $$\p-\z_2\t (\u_2-\u_1)-(\z_2-\z_1)\t \u_1\p_{H^{-1/2}(\O)}\le \bar{S}(\p z_2\p_{L^2(\O)}\p\u_2-\u_1\p_{L^4(\O)}+\p z_2-z_1\p_{L^2(\O)}\p \u_1\p_{L^4(\O)}).$$
   Finally substituting this last inequality and the bound (\ref{5bliu2u1}) in (\ref{5dC3/2}), using (\ref{5dSinf}) and (\ref{1dSp}), give (\ref{5bliu2u1}).\hfill$\diamondsuit$\\[0.2cm]
  \H We now can prove the two previous theorems. Let $(\u_1,p_1,z_1)$ and $(\u_2,p_2,z_2)$ be two solutions of Problem (\ref{2eqe}). For $i=1,2$, $(\u_i,z_i)$ is  a solution of the transport equation (\ref{1et}). So we have, for $i=1,2$,
  $$ \nu z_i +\a \u_i\,.\,\n z_i=\nu\,\cs\u_i+\a\,\cs\f.$$
  Then $z_2-z_1$ is a solution of the following equation
  $$\nu(z_2-z_1)+\a\u_2\,.\,\n(z_2-z_1)=\nu\,\cs(\u_2-\u_1)-\a(\u_2-\u_1)\,.\,\n z_1.$$
  Now, taking the scalar product of both sides of the previous equation with $z_2-z_1$, we obtain:
  $$\nu\p z_2-z_1\p_{L^2(\O)}^2+\a(\u_2\,.\,\n(z_2-z_1),z_2-z_1)=\nu\,(\cs(\u_2-\u_1),z_2-z_1)-\a((\u_2-\u_1)\,.\,\n z_1,z_2-z_1).$$
  \H Let us note that $z_2-z_1$ belongs to $X_{\u_2}(\G^-)$ (see (\ref{1dg-}), (\ref{1dXug0})), therefore, owing to the Proposition \ref{2pgreenfg-}, we have $\a(\u_2\,.\,\n(z_2-z_1),z_2-z_1)\ge 0$. Hence, we derive
  \begin{equation}\label{5eiz1z2}
  \nu\p z_2-z_1\p_{L^2(\O)}\le \nu\sqrt{2}|\u_2-\u_1|_{H^1(\O)}+|\a|\p \u_2-\u_1\p_{L^{\iy}(\O)}\p\n z_1\p_{L^2(\O)}.\end{equation}
  Next, in view of (\ref{5bH1u2u1}) and (\ref{5bliu2u1}), we obtain the following estimate
  \begin{equation}
  \p z_2-z_1\p_{L^2(\O)}\le \F{S^*_4}{\nu}\p\u_1\p_{H^1(\O)}(\sqrt{2}\,S_4
  +\F{|\a|}{\nu}C_{3/2} S_{\iy} \bar{S}(1+\F{S_4^2}{\nu}\p z_2\p_{L^2(\O)})\p\n z_1\p_{L^2(\O)})\p z_2-z_1\p_{L^2(\O)}.
  \label{5iiuf}\end{equation}
  Setting
  \begin{equation}\label{5dL1}
  L_1(\a,\nu,\f,\g,h)=\F{2|\a|}{\nu}\p\cs\f\p_{L^2(\O)}
   +4\p z_h \p_{L^2(\O)}+C(\a,\nu,\f,\g,h),\end{equation}
   and
 \begin{equation}\label{5dL2}
  L_2(\a,\nu,\f,\g,h)=\F{S_2\sqrt{S_2^2+1}}{\nu}\p\f\p_{L^2(\O)}
   +T\p\g\p_{H^{1/2}(\P\O)}\left(1+
   \F{S_4S^*_4\sqrt{S_2^2+1}}{\nu}\,
   L_1(\a,\nu,\f,\g,\h)\right),\end{equation}
   where the function $C$ is defined by (\ref{3dC}), the bounds (\ref{2mz}) and (\ref{2mu}) give
   \begin{equation}\label{5mzL2uH1}
  \p z_2\p_{L^2(\O)}\le  L_1(\a,\nu,\f,\g,h),\H\p\u_1\p_{H^1(\O)}\le L_2(\a,\nu,\f,\g,h).\end{equation}
   With these notations, from (\ref{5iiuf}), we derive
   \begin{eqnarray}
    \p z_2-z_1\p_{L^2(\O)}(1- \F{S^*_4}{\nu}L_2(\a,\nu,\f,\g,h)(\sqrt{2}\,S_4\nonumber\\
  +\F{|\a|}{\nu}C_{3/2} S_{\iy} \bar{S}(1+\F{S_4^2}{\nu}L_1(\a,\nu,\f,\g,h))\p\n z_1\p_{L^2(\O)}))\le 0.\label{5iufin}
  \end{eqnarray}
  Hence, in view of (\ref{5cuz}), we obtain $z_1=z_2$, which implies that Problem \ref{2eqe} has a unique solution and which ends the proof of Theorem \ref{5thupb2}.\\
  \H Finally, applying Lemma \ref{5lbzH1bis} yields
  $$\p\n z_1\p_{L^2(\O)}\le L_3(\a,\nu,\f,\g,h),$$
  which gives with (\ref{5iufin})
  \begin{eqnarray*}
    \p z_2-z_1\p_{L^2(\O)}(1- \F{S^*_4}{\nu}L_2(\a,\nu,\f,\g,h)(\sqrt{2}\,S_4\\
  +\F{|\a|}{\nu}C_{3/2} S_{\iy} \bar{S}(1+\F{S_4^2}{\nu}L_1(\a,\nu,\f,\g,h))\,L_3(\a,\nu,\f,\g,h)))\le 0.
  \end{eqnarray*}
  Then, in view of (\ref{5cuzb}), we  again derive $z_1=z_2$, which implies uniqueness and ends the proof of Theorem \ref{5thupb2b}.\hfill$\diamondsuit$\\[0.5cm]
  \textbf{References}\\[0.2cm]
  \lk 1\rk\hs D.N. Arnold, L.R. Scott and M. Vogelius, "Regular inversion of
  the divergence operator with Dirichlet boundary conditions on a polygon",
  \textit{Ann. della Scuola Norm. sup. di Pisa}, Serie IV, XV, Fasc. II,
  169-192 (1988).\\[0.2cm]
  \lk 2\rk\hs J.M. Bernard, \textit{Fluides de second et troisi\`eme grade en 
  dimension trois: solution globale et r\'egularit\'e, Th\`ese de doctorat,} 
  Universit\'e Pierre-et-Marie-Curie (1998). \\[0.2cm]
  \lk 3\rk\hs J.M. Bernard, "Stationary Problem of Second-grade Fluids in
  Three Dimensions: \linebreak Existence, Uniqueness and Regularity", Math. Meth.
  Appl. Sci., \tb{22}, 655-687 (1999).\\[0.2cm]
  \lk 4\rk\hs J.M. Bernard, "Solutions $W^{2,p},\ p>3$, of equations of a family of second grade
  fluids with a boundary of class $C^{1,1}$", Communications on Applied
  Nonlinear Analysis. Vol. 9, No. 1, 1-29 (2002).\\[0.2cm]
  \lk 5\rk\hs J.M. Bernard, "Steady transport equations in the case where the normal component of the 
  velocity does not vanish on the boundary", SIAM J. Math. Anal. Vol. 44, No. 2, 993-1018 (2012).\\[0.2cm]
  \lk 6\rk\hs J.M. Bernard, "Problem of second grade fluids in convex polyhedrons",
  SIAM J. Math. Anal. Vol. 44, No. 3, 2018-2038 (2012).\\[0.2cm]
  \lk 7\rk\hs J.M. Bernard, "Solutions $H^1$ of the steady transport equation in the case where the normal component 
  of the velocity does not vanish on    the boundary'', 
   J. Math. Pures  Appl. \tb{78}, 10, 981-1011 (2016).\\[0.2cm]
  \l 8\r\hs D. Bresch and J. Lemoine, "Stationary solutions for second-grade 
  fluids equations", $M.^3 A.S.$, \textbf{8}, $\mr{n}^o 5$ (1998).\\[0.2cm] 
  \lk 9\rk\hs D. Bresch and J. Lemoine, "On the existence of solutions  
  for non-stationary third-grade fluids", Int. J. Non-linear Mechanics, 
  \textbf{34}, 3, 485-498 (1998).\\[0.2cm] 
  \lk 10\rk\hs D. Cioranescu and V. Girault, "Weak and classical solutions of a family 
  of second grade fluids", Int. J. Non-linear Mechanics, \textbf{32}, 2, 
  317-335 (1997).\\[0.2cm]  
  \lk 11\rk\hs D. Cioranescu and E.H. Ouazar, "Existence et unicit\'e pour les fluides
  de second grade", Note CRAS \textbf{298} S\'erie I, 285-287 (1984).\\[0.2cm] 
  \lk 12\rk\hs D. Cioranescu and E.H. Ouazar, "Existence and uniqueness for fluids of 
  second grade", in Nonlinear Partial Differential Equations, Coll\`ege de 
  France Seminar, Pitman \textbf{109} 178-197 (1984).\\[0.2cm]  
  \lk 13\rk\hs V. Coscia and G.P. Galdi, "Existence, uniqueness and stability of regular 
  steady motions of second-grade fluid", Int. J. Non-Linear Mech. 
  \textbf{29}(4), 493-506 (1994).\\[0.2cm]
  \lk 14\rk\hs M. Dauge, "Stationary Stokes and Navier-Stokes systems on
  two- or three-dimensional domains with corners. Part I: linearized
  equations", SIAM J. MATH. Anal. Vol. 20, No. 1, 74-97 (1989).\\[0.2cm]
  \lk 15\rk\hs R.J. Diperna and P.L. Lions, "Ordinary differential equations,
  transport theory and Sobolev spaces", Invent. Math. 98, 511-547 (1989).\\[0.2cm] 
  \lk 16\rk\hs J.E. Dunn and R.L. Fosdick, "Thermodynamics, stability and boundedness of 
  fluids of complexity two and fluids of second grade", Arch. Rat. Mech. Anal.
  \textbf{56},3, 191-252 (1974).\\[0.2cm]
  \lk 17\rk\hs J.E. Dunn and K.R. Rajagopal, "Fluids of differential type:
  Critical review and thermodynamic analysis", Int. J. Engng. Sci \tb{33},
  5, 689-729 (1995).\\[0.2cm] 
  \lk 18\rk\hs R.L. Fosdick and K.R. Rajagopal, "Anomalous features in the model of 
  second order fluids", Arch. Rat. Mech. Anal. \textbf{70}, 3, 1-46 (1979).
  \\[0.2cm]  
  \lk 19\rk\hs R.L. Fosdick and K.R. Rajagopal, "Thermodynamics and stability of fluids 
  of third grade", Proc. Royal Soc. London \textbf{A 339}, 351-377 (1980).
  \\[0.2cm] 
  \lk 20\rk\hs G.P. Galdi, M. Grobbelaar-Van Dalsen and N. Sauer, "Existence and 
  uniqueness of classical solutions of the equations of motion for second 
  grade fluids", Arch. Rat. Mech. Anal. \textbf{124}, 221-237 (1993).\\[0.2cm] 
  \lk 21\rk\hs G.P. Galdi and A. Sequeira, "Further existence results for classical 
  solutions of the equations of second grade fluids", Arch. Rat. Mech. Anal. 
  \textbf{128}, 297-312 (1994).\\[0.2cm]
  \lk 22\rk\hs V. Girault and P.A. Raviart, \textit{Finite Element Methods for
  the Navier-Stokes Equations. Theory and Algorithms}, SCM 5, Spinger-Verlag,
  Berlin, 1986.\\[0.2cm]
  \lk 23\rk\hs V. Girault and L.R. Scott, "Analysis of a two-dimensional
  grade-two fluid model with a tangential boundary condition", J. Math. Pures
  Appl. \tb{78}, 10, 981-1011 (1999).\\[0.2cm]
  \lk 24\rk\hs V. Girault and L.R. Scott, "Finite-element discretizations of a
  two-dimensional grade-two fluid model", M$^2$AN, \tb{35},
  1007-1053 (2001).\\[0.2cm]
  \lk 25\rk\hs P. Grisvard, \textit{Elliptic Problems in Nonsmooth Domains},
  Pitman Monographs and Stu-\linebreak dies in Mathematics 24, Pitman, Boston, MA,
  1985.\\[0.2cm]
  \lk 26\rk\hs J. Ne$\check{\mathrm{c}}$as, \textit{Les M\'ethodes Directes en Th\'eorie des 
  Equations Elliptiques}, Masson, Paris (1967).\\[0.2cm] 
  \lk 27\rk\hs W. Noll and C. Truesdell, \textit{The Nonlinear Field Theory of 
  Mechanics. Handbuch of Physik,} Vol. III, Springer-Verlag, Berlin (1975).
  \\[0.2cm] 
  \lk 28\rk\hs E.H. Ouazar, \textit{Sur les Fluides de Second Grade, Th\`ese 3\`eme 
  Cycle,} Universit\'e Pierre-et-Marie-Curie (1981).\\[0.2cm]
  \lk 29\rk\hs J.H. Videman, \textit{Mathematical analysis of viscoelastic 
  non-newtonian fluids, Thesis}, Universit\'e de Lisbonne (1997).

\end{document}